\documentclass{gtart}

\def\ifplaintex{\expandafter\ifx\csname documentclass\endcsname\relax}

\def\gtp{{\mathsurround=0pt\it $\cal G\mskip-2mu$eometry \&\ 
$\cal T\!\!$opology $\cal P\!$ublications}}  

\def\recd{{\small Received:\qua\receiveddate\ifx\reviseddate\relax
\else\qquad Revised:\qua\reviseddate\fi\par}} 

\def\lognumber#1{\def\thelognumber{#1}}
\def\volumenumber#1{\def\thevolumenumber{#1}}
\def\volumeyear#1{\def\thevolumeyear{#1}}
\def\papernumber#1{\def\thepapernumber{#1}}
\def\pagenumbers#1#2{\def\startpage{#1}\def\finishpage{#2}}
\def\published#1{\def\publishdate{#1}}

\def\received#1{\def\receiveddate{#1}}

\def\accepted#1{\def\accepteddate{#1}}
\def\asciititle#1{\def\theasciititle{#1}}

\long\def\asciiabstract#1{\long\def\theasciiabstract{#1}}

\let\\\par\let\thelognumber\relax\let\thevolumenumber\relax
\let\thepapernumber\relax\let\thevolumeyear\relax\let\startpage\relax
\let\finishpage\relax\let\publishdate\relax\let\receiveddate\relax
\let\reviseddate\relax\let\accepteddate\relax\let\theasciititle\relax
\let\theasciiauthors\relax
\let\theasciiabstract\relax

\let\theasciiemail\relax

\ifplaintex
\font\logobig=cmssbx10 scaled 3836
\font\logomed=cmssbx10 scaled 2557
\else
\font\logobig=cmssbx10 scaled 4200
\font\logomed=cmssbx10 scaled 2800
\fi

\long\def\makeagttitle{   
\count0=\startpage
\agt\hfill      
\hbox to 45truept{\vbox to 0pt{\vglue -13truept{\logomed A\kern -.37em{\logobig 
T}\kern -.38em G}\vss}\hss}
\break
{\small Volume \thevolumenumber\ (\thevolumeyear)
\startpage--\finishpage\nl
Published: \publishdate}

\vglue .25truein

{\parskip=0pt\leftskip 0pt plus
1fil\def\\{\par\smallskip}{\Large\bf\thetitle}\par\medskip} \vglue
0.05truein

{\parskip=0pt\leftskip 0pt plus 1fil\def\\{\par}{\sc\theauthors}
\par\medskip}%
 
\vglue 0.03truein 

{\small\leftskip 25truept\rightskip 25truept{\bf Abstract}\stdspace\theabstract

{\bf AMS Classification}\stdspace\theprimaryclass
\ifx\thesecondaryclass\relax\else; \thesecondaryclass\fi\par
{\bf Keywords}\stdspace \thekeywords\par}\vglue 7truept

}   

\ifplaintex
\font\phead=cmsl9 scaled 950
\font\pnum=cmbx10 scaled 913
\font\pfoot=cmsl9 scaled 950
\headline{\vbox to 0pt{\vskip -4.5mm\line{\small\phead\ifnum
\count0=\startpage ISSN 1472-2739 (on-line) 1472-2747 (printed)
\hfill {\pnum\folio}\else\ifodd\count0\def\\{ }%
\ifx\theshorttitle\relax\thetitle\else\theshorttitle\fi\hfill{\pnum\folio}
\else\def\\{ and }{\pnum\folio}\hfill\ifx\theshortauthors\relax\theauthors
\else\theshortauthors\fi\fi\fi}\vss}}
\footline{\vbox to 0pt{\vglue 0mm\line{\small\pfoot\ifnum\count0=\startpage
\copyright\ \gtp\hfill\else
\agt, Volume \thevolumenumber\ (\thevolumeyear)\hfill\fi}\vss}}
\else
\font=cmsl9 scaled 1050
\font\lnum=cmbx10 
\font=cmsl9 scaled 1050
\makeatletter
\def\@oddhead{{\small\lhead\ifnum\count0=\startpage ISSN 1472-2739 
(on-line) 1472-2747 (printed)\hfill {\lnum\number\count0}\else\ifodd\count0
\def\\{ }\ifx\theshorttitle\relax \thetitle \else\theshorttitle\fi\hfill
{\lnum\number\count0}\else\def\\{ and }{\lnum\number\count0}
\hfill\ifx\theshortauthors\relax 
\theauthors\else\theshortauthors\fi\fi\fi}}\def\@evenhead{\@oddhead}
\def\@oddfoot{\small\lfoot\ifnum\count0=\startpage\copyright\ \gtp\hfill\else
\agt, Volume \thevolumenumber\ (\thevolumeyear)\hfill\fi}
\def\@evenfoot{\@oddfoot}
\makeatother
\fi
\let\maketitlepage\makeagttitle
\let\makeshorttitle\maketitlepage
\let\maketitle\maketitlepage

\newwrite\gtoutfile
\long\gdef\makeheadfile{  
{\def\\{, }\def\s{ }
\immediate\openout\gtoutfile head.xxx
\immediate\write\gtoutfile{To: math@arxiv.org}
\immediate\write\gtoutfile{Subject: put OR rep NNNNN:ppppp}
\immediate\write\gtoutfile{--text follows this line--}
\immediate\write\gtoutfile{Proxy-for: \ifx\theasciiauthors\relax
\theauthors\else\theasciiauthors\fi\s<\ifx\theasciiemail\relax\theemail\else\theasciiemail\fi>}
\immediate\write\gtoutfile{\noexpand\\}
\immediate\write\gtoutfile{Authors: \ifx\theasciiauthors\relax
\theauthors\else\theasciiauthors\fi}
{\def\\{ }\immediate\write\gtoutfile{Title: \ifx\theasciititle\relax
\thetitle\else\theasciititle\fi}}
\immediate\write\gtoutfile{Subj-class: GT or SG, GR etc}
\immediate\write\gtoutfile{MSC-class: \theprimaryclass\ifx\thesecondaryclass\relax\else, \thesecondaryclass\fi}
\immediate\write\gtoutfile{Journal-ref: Algebr. Geom. Topol. \thevolumenumber\s
(\thevolumeyear) \startpage-\finishpage}
\immediate\write\gtoutfile{Comments: Published by Algebraic and
Geometric Topology at}
\immediate\write\gtoutfile{\s\s\s  http://www.maths.warwick.ac.uk/agt/AGTVol\thevolumenumber/agt-\thevolumenumber-\thepapernumber.abs.html}
\immediate\write\gtoutfile{\noexpand\\}
\immediate\write\gtoutfile{}
\ifx\theasciiabstract\relax
\immediate\write\gtoutfile{\theabstract}\else
\immediate\write\gtoutfile{\theasciiabstract}\fi
\immediate\write\gtoutfile{}
\immediate\write\gtoutfile{\noexpand\\}
\immediate\write\gtoutfile{}
\immediate\closeout\gtoutfile}}  

\def\maketitlepage{\makeagttitle\makeheadfile}
\let\makeshorttitle\maketitlepage
\let\maketitle\maketitlepage

\def\ifplaintex{\expandafter\ifx\csname documentclass\endcsname\relax}

\def\gtp{{\mathsurround=0pt\it $\cal G\mskip-2mu$eometry \&\ 
$\cal T\!\!$opology $\cal P\!$ublications}}  

\def\recd{{\small Received:\qua\receiveddate\ifx\reviseddate\relax
\else\qquad Revised:\qua\reviseddate\fi\par}} 

\def\lognumber#1{\def\thelognumber{#1}}
\def\volumenumber#1{\def\thevolumenumber{#1}}
\def\volumeyear#1{\def\thevolumeyear{#1}}
\def\papernumber#1{\def\thepapernumber{#1}}
\def\pagenumbers#1#2{\def\startpage{#1}\def\finishpage{#2}}
\def\published#1{\def\publishdate{#1}}

\def\received#1{\def\receiveddate{#1}}

\def\accepted#1{\def\accepteddate{#1}}
\def\asciititle#1{\def\theasciititle{#1}}

\long\def\asciiabstract#1{\long\def\theasciiabstract{#1}}

\let\\\par\let\thelognumber\relax\let\thevolumenumber\relax
\let\thepapernumber\relax\let\thevolumeyear\relax\let\startpage\relax
\let\finishpage\relax\let\publishdate\relax\let\receiveddate\relax
\let\reviseddate\relax\let\accepteddate\relax\let\theasciititle\relax
\let\theasciiauthors\relax
\let\theasciiabstract\relax

\let\theasciiemail\relax

\ifplaintex
\font\logobig=cmssbx10 scaled 3836
\font\logomed=cmssbx10 scaled 2557
\else
\font\logobig=cmssbx10 scaled 4200
\font\logomed=cmssbx10 scaled 2800
\fi

\long\def\makeagttitle{   
\count0=\startpage
\agt\hfill      
\hbox to 45truept{\vbox to 0pt{\vglue -13truept{\logomed A\kern -.37em{\logobig 
T}\kern -.38em G}\vss}\hss}
\break
{\small Volume \thevolumenumber\ (\thevolumeyear)
\startpage--\finishpage\nl
Published: \publishdate}

\vglue .25truein

{\parskip=0pt\leftskip 0pt plus
1fil\def\\{\par\smallskip}{\Large\bf\thetitle}\par\medskip} \vglue
0.05truein

{\parskip=0pt\leftskip 0pt plus 1fil\def\\{\par}{\sc\theauthors}
\par\medskip}%
 
\vglue 0.03truein 

{\small\leftskip 25truept\rightskip 25truept{\bf Abstract}\stdspace\theabstract

{\bf AMS Classification}\stdspace\theprimaryclass
\ifx\thesecondaryclass\relax\else; \thesecondaryclass\fi\par
{\bf Keywords}\stdspace \thekeywords\par}\vglue 7truept

}   

\ifplaintex
\font\phead=cmsl9 scaled 950
\font\pnum=cmbx10 scaled 913
\font\pfoot=cmsl9 scaled 950
\headline{\vbox to 0pt{\vskip -4.5mm\line{\small\phead\ifnum
\count0=\startpage ISSN 1472-2739 (on-line) 1472-2747 (printed)
\hfill {\pnum\folio}\else\ifodd\count0\def\\{ }%
\ifx\theshorttitle\relax\thetitle\else\theshorttitle\fi\hfill{\pnum\folio}
\else\def\\{ and }{\pnum\folio}\hfill\ifx\theshortauthors\relax\theauthors
\else\theshortauthors\fi\fi\fi}\vss}}
\footline{\vbox to 0pt{\vglue 0mm\line{\small\pfoot\ifnum\count0=\startpage
\copyright\ \gtp\hfill\else
\agt, Volume \thevolumenumber\ (\thevolumeyear)\hfill\fi}\vss}}
\else
\font=cmsl9 scaled 1050
\font\lnum=cmbx10 
\font=cmsl9 scaled 1050
\makeatletter
\def\@oddhead{{\small\lhead\ifnum\count0=\startpage ISSN 1472-2739 
(on-line) 1472-2747 (printed)\hfill {\lnum\number\count0}\else\ifodd\count0
\def\\{ }\ifx\theshorttitle\relax \thetitle \else\theshorttitle\fi\hfill
{\lnum\number\count0}\else\def\\{ and }{\lnum\number\count0}
\hfill\ifx\theshortauthors\relax 
\theauthors\else\theshortauthors\fi\fi\fi}}\def\@evenhead{\@oddhead}
\def\@oddfoot{\small\lfoot\ifnum\count0=\startpage\copyright\ \gtp\hfill\else
\agt, Volume \thevolumenumber\ (\thevolumeyear)\hfill\fi}
\def\@evenfoot{\@oddfoot}
\makeatother
\fi
\let\maketitlepage\makeagttitle
\let\makeshorttitle\maketitlepage
\let\maketitle\maketitlepage

\newwrite\gtoutfile
\long\gdef\makeheadfile{  
{\def\\{, }\def\s{ }
\immediate\openout\gtoutfile head.xxx
\immediate\write\gtoutfile{To: math@arxiv.org}
\immediate\write\gtoutfile{Subject: put OR rep NNNNN:ppppp}
\immediate\write\gtoutfile{--text follows this line--}
\immediate\write\gtoutfile{Proxy-for: \ifx\theasciiauthors\relax
\theauthors\else\theasciiauthors\fi\s<\ifx\theasciiemail\relax\theemail\else\theasciiemail\fi>}
\immediate\write\gtoutfile{\noexpand\\}
\immediate\write\gtoutfile{Authors: \ifx\theasciiauthors\relax
\theauthors\else\theasciiauthors\fi}
{\def\\{ }\immediate\write\gtoutfile{Title: \ifx\theasciititle\relax
\thetitle\else\theasciititle\fi}}
\immediate\write\gtoutfile{Subj-class: GT or SG, GR etc}
\immediate\write\gtoutfile{MSC-class: \theprimaryclass\ifx\thesecondaryclass\relax\else, \thesecondaryclass\fi}
\immediate\write\gtoutfile{Journal-ref: Algebr. Geom. Topol. \thevolumenumber\s
(\thevolumeyear) \startpage-\finishpage}
\immediate\write\gtoutfile{Comments: Published by Algebraic and
Geometric Topology at}
\immediate\write\gtoutfile{\s\s\s  http://www.maths.warwick.ac.uk/agt/AGTVol\thevolumenumber/agt-\thevolumenumber-\thepapernumber.abs.html}
\immediate\write\gtoutfile{\noexpand\\}
\immediate\write\gtoutfile{}
\ifx\theasciiabstract\relax
\immediate\write\gtoutfile{\theabstract}\else
\immediate\write\gtoutfile{\theasciiabstract}\fi
\immediate\write\gtoutfile{}
\immediate\write\gtoutfile{\noexpand\\}
\immediate\write\gtoutfile{}
\immediate\closeout\gtoutfile}}  

\def\maketitlepage{\makeagttitle\makeheadfile}
\let\makeshorttitle\maketitlepage
\let\maketitle\maketitlepage

\def\ifplaintex{\expandafter\ifx\csname documentclass\endcsname\relax}

\def\gtp{{\mathsurround=0pt\it $\cal G\mskip-2mu$eometry \&\ 
$\cal T\!\!$opology $\cal P\!$ublications}}  

\def\recd{{\small Received:\qua\receiveddate\ifx\reviseddate\relax
\else\qquad Revised:\qua\reviseddate\fi\par}} 

\def\lognumber#1{\def\thelognumber{#1}}
\def\volumenumber#1{\def\thevolumenumber{#1}}
\def\volumeyear#1{\def\thevolumeyear{#1}}
\def\papernumber#1{\def\thepapernumber{#1}}
\def\pagenumbers#1#2{\def\startpage{#1}\def\finishpage{#2}}
\def\published#1{\def\publishdate{#1}}

\def\received#1{\def\receiveddate{#1}}

\def\accepted#1{\def\accepteddate{#1}}
\def\asciititle#1{\def\theasciititle{#1}}

\long\def\asciiabstract#1{\long\def\theasciiabstract{#1}}

\let\\\par\let\thelognumber\relax\let\thevolumenumber\relax
\let\thepapernumber\relax\let\thevolumeyear\relax\let\startpage\relax
\let\finishpage\relax\let\publishdate\relax\let\receiveddate\relax
\let\reviseddate\relax\let\accepteddate\relax\let\theasciititle\relax
\let\theasciiauthors\relax
\let\theasciiabstract\relax

\let\theasciiemail\relax

\ifplaintex
\font\logobig=cmssbx10 scaled 3836
\font\logomed=cmssbx10 scaled 2557
\else
\font\logobig=cmssbx10 scaled 4200
\font\logomed=cmssbx10 scaled 2800
\fi

\long\def\makeagttitle{   
\count0=\startpage
\agt\hfill      
\hbox to 45truept{\vbox to 0pt{\vglue -13truept{\logomed A\kern -.37em{\logobig 
T}\kern -.38em G}\vss}\hss}
\break
{\small Volume \thevolumenumber\ (\thevolumeyear)
\startpage--\finishpage\nl
Published: \publishdate}

\vglue .25truein

{\parskip=0pt\leftskip 0pt plus
1fil\def\\{\par\smallskip}{\Large\bf\thetitle}\par\medskip} \vglue
0.05truein

{\parskip=0pt\leftskip 0pt plus 1fil\def\\{\par}{\sc\theauthors}
\par\medskip}%
 
\vglue 0.03truein 

{\small\leftskip 25truept\rightskip 25truept{\bf Abstract}\stdspace\theabstract

{\bf AMS Classification}\stdspace\theprimaryclass
\ifx\thesecondaryclass\relax\else; \thesecondaryclass\fi\par
{\bf Keywords}\stdspace \thekeywords\par}\vglue 7truept

}   

\ifplaintex
\font\phead=cmsl9 scaled 950
\font\pnum=cmbx10 scaled 913
\font\pfoot=cmsl9 scaled 950
\headline{\vbox to 0pt{\vskip -4.5mm\line{\small\phead\ifnum
\count0=\startpage ISSN 1472-2739 (on-line) 1472-2747 (printed)
\hfill {\pnum\folio}\else\ifodd\count0\def\\{ }%
\ifx\theshorttitle\relax\thetitle\else\theshorttitle\fi\hfill{\pnum\folio}
\else\def\\{ and }{\pnum\folio}\hfill\ifx\theshortauthors\relax\theauthors
\else\theshortauthors\fi\fi\fi}\vss}}
\footline{\vbox to 0pt{\vglue 0mm\line{\small\pfoot\ifnum\count0=\startpage
\copyright\ \gtp\hfill\else
\agt, Volume \thevolumenumber\ (\thevolumeyear)\hfill\fi}\vss}}
\else
\font=cmsl9 scaled 1050
\font\lnum=cmbx10 
\font=cmsl9 scaled 1050
\makeatletter
\def\@oddhead{{\small\lhead\ifnum\count0=\startpage ISSN 1472-2739 
(on-line) 1472-2747 (printed)\hfill {\lnum\number\count0}\else\ifodd\count0
\def\\{ }\ifx\theshorttitle\relax \thetitle \else\theshorttitle\fi\hfill
{\lnum\number\count0}\else\def\\{ and }{\lnum\number\count0}
\hfill\ifx\theshortauthors\relax 
\theauthors\else\theshortauthors\fi\fi\fi}}\def\@evenhead{\@oddhead}
\def\@oddfoot{\small\lfoot\ifnum\count0=\startpage\copyright\ \gtp\hfill\else
\agt, Volume \thevolumenumber\ (\thevolumeyear)\hfill\fi}
\def\@evenfoot{\@oddfoot}
\makeatother
\fi
\let\maketitlepage\makeagttitle
\let\makeshorttitle\maketitlepage
\let\maketitle\maketitlepage

\newwrite\gtoutfile
\long\gdef\makeheadfile{  
{\def\\{, }\def\s{ }
\immediate\openout\gtoutfile head.xxx
\immediate\write\gtoutfile{To: math@arxiv.org}
\immediate\write\gtoutfile{Subject: put OR rep NNNNN:ppppp}
\immediate\write\gtoutfile{--text follows this line--}
\immediate\write\gtoutfile{Proxy-for: \ifx\theasciiauthors\relax
\theauthors\else\theasciiauthors\fi\s<\ifx\theasciiemail\relax\theemail\else\theasciiemail\fi>}
\immediate\write\gtoutfile{\noexpand\\}
\immediate\write\gtoutfile{Authors: \ifx\theasciiauthors\relax
\theauthors\else\theasciiauthors\fi}
{\def\\{ }\immediate\write\gtoutfile{Title: \ifx\theasciititle\relax
\thetitle\else\theasciititle\fi}}
\immediate\write\gtoutfile{Subj-class: GT or SG, GR etc}
\immediate\write\gtoutfile{MSC-class: \theprimaryclass\ifx\thesecondaryclass\relax\else, \thesecondaryclass\fi}
\immediate\write\gtoutfile{Journal-ref: Algebr. Geom. Topol. \thevolumenumber\s
(\thevolumeyear) \startpage-\finishpage}
\immediate\write\gtoutfile{Comments: Published by Algebraic and
Geometric Topology at}
\immediate\write\gtoutfile{\s\s\s  http://www.maths.warwick.ac.uk/agt/AGTVol\thevolumenumber/agt-\thevolumenumber-\thepapernumber.abs.html}
\immediate\write\gtoutfile{\noexpand\\}
\immediate\write\gtoutfile{}
\ifx\theasciiabstract\relax
\immediate\write\gtoutfile{\theabstract}\else
\immediate\write\gtoutfile{\theasciiabstract}\fi
\immediate\write\gtoutfile{}
\immediate\write\gtoutfile{\noexpand\\}
\immediate\write\gtoutfile{}
\immediate\closeout\gtoutfile}}  

\def\maketitlepage{\makeagttitle\makeheadfile}
\let\makeshorttitle\maketitlepage
\let\maketitle\maketitlepage

\def\ifplaintex{\expandafter\ifx\csname documentclass\endcsname\relax}

\def\gtp{{\mathsurround=0pt\it $\cal G\mskip-2mu$eometry \&\ 
$\cal T\!\!$opology $\cal P\!$ublications}}  

\def\recd{{\small Received:\qua\receiveddate\ifx\reviseddate\relax
\else\qquad Revised:\qua\reviseddate\fi\par}} 

\def\lognumber#1{\def\thelognumber{#1}}
\def\volumenumber#1{\def\thevolumenumber{#1}}
\def\volumeyear#1{\def\thevolumeyear{#1}}
\def\papernumber#1{\def\thepapernumber{#1}}
\def\pagenumbers#1#2{\def\startpage{#1}\def\finishpage{#2}}
\def\published#1{\def\publishdate{#1}}

\def\received#1{\def\receiveddate{#1}}

\def\accepted#1{\def\accepteddate{#1}}
\def\asciititle#1{\def\theasciititle{#1}}

\long\def\asciiabstract#1{\long\def\theasciiabstract{#1}}

\let\\\par\let\thelognumber\relax\let\thevolumenumber\relax
\let\thepapernumber\relax\let\thevolumeyear\relax\let\startpage\relax
\let\finishpage\relax\let\publishdate\relax\let\receiveddate\relax
\let\reviseddate\relax\let\accepteddate\relax\let\theasciititle\relax
\let\theasciiauthors\relax
\let\theasciiabstract\relax

\let\theasciiemail\relax

\ifplaintex
\font\logobig=cmssbx10 scaled 3836
\font\logomed=cmssbx10 scaled 2557
\else
\font\logobig=cmssbx10 scaled 4200
\font\logomed=cmssbx10 scaled 2800
\fi

\long\def\makeagttitle{   
\count0=\startpage
\agt\hfill      
\hbox to 45truept{\vbox to 0pt{\vglue -13truept{\logomed A\kern -.37em{\logobig 
T}\kern -.38em G}\vss}\hss}
\break
{\small Volume \thevolumenumber\ (\thevolumeyear)
\startpage--\finishpage\nl
Published: \publishdate}

\vglue .25truein

{\parskip=0pt\leftskip 0pt plus
1fil\def\\{\par\smallskip}{\Large\bf\thetitle}\par\medskip} \vglue
0.05truein

{\parskip=0pt\leftskip 0pt plus 1fil\def\\{\par}{\sc\theauthors}
\par\medskip}%
 
\vglue 0.03truein 

{\small\leftskip 25truept\rightskip 25truept{\bf Abstract}\stdspace\theabstract

{\bf AMS Classification}\stdspace\theprimaryclass
\ifx\thesecondaryclass\relax\else; \thesecondaryclass\fi\par
{\bf Keywords}\stdspace \thekeywords\par}\vglue 7truept

}   

\ifplaintex
\font\phead=cmsl9 scaled 950
\font\pnum=cmbx10 scaled 913
\font\pfoot=cmsl9 scaled 950
\headline{\vbox to 0pt{\vskip -4.5mm\line{\small\phead\ifnum
\count0=\startpage ISSN 1472-2739 (on-line) 1472-2747 (printed)
\hfill {\pnum\folio}\else\ifodd\count0\def\\{ }%
\ifx\theshorttitle\relax\thetitle\else\theshorttitle\fi\hfill{\pnum\folio}
\else\def\\{ and }{\pnum\folio}\hfill\ifx\theshortauthors\relax\theauthors
\else\theshortauthors\fi\fi\fi}\vss}}
\footline{\vbox to 0pt{\vglue 0mm\line{\small\pfoot\ifnum\count0=\startpage
\copyright\ \gtp\hfill\else
\agt, Volume \thevolumenumber\ (\thevolumeyear)\hfill\fi}\vss}}
\else
\font=cmsl9 scaled 1050
\font\lnum=cmbx10 
\font=cmsl9 scaled 1050
\makeatletter
\def\@oddhead{{\small\lhead\ifnum\count0=\startpage ISSN 1472-2739 
(on-line) 1472-2747 (printed)\hfill {\lnum\number\count0}\else\ifodd\count0
\def\\{ }\ifx\theshorttitle\relax \thetitle \else\theshorttitle\fi\hfill
{\lnum\number\count0}\else\def\\{ and }{\lnum\number\count0}
\hfill\ifx\theshortauthors\relax 
\theauthors\else\theshortauthors\fi\fi\fi}}\def\@evenhead{\@oddhead}
\def\@oddfoot{\small\lfoot\ifnum\count0=\startpage\copyright\ \gtp\hfill\else
\agt, Volume \thevolumenumber\ (\thevolumeyear)\hfill\fi}
\def\@evenfoot{\@oddfoot}
\makeatother
\fi
\let\maketitlepage\makeagttitle
\let\makeshorttitle\maketitlepage
\let\maketitle\maketitlepage

\newwrite\gtoutfile
\long\gdef\makeheadfile{  
{\def\\{, }\def\s{ }
\immediate\openout\gtoutfile head.xxx
\immediate\write\gtoutfile{To: math@arxiv.org}
\immediate\write\gtoutfile{Subject: put OR rep NNNNN:ppppp}
\immediate\write\gtoutfile{--text follows this line--}
\immediate\write\gtoutfile{Proxy-for: \ifx\theasciiauthors\relax
\theauthors\else\theasciiauthors\fi\s<\ifx\theasciiemail\relax\theemail\else\theasciiemail\fi>}
\immediate\write\gtoutfile{\noexpand\\}
\immediate\write\gtoutfile{Authors: \ifx\theasciiauthors\relax
\theauthors\else\theasciiauthors\fi}
{\def\\{ }\immediate\write\gtoutfile{Title: \ifx\theasciititle\relax
\thetitle\else\theasciititle\fi}}
\immediate\write\gtoutfile{Subj-class: GT or SG, GR etc}
\immediate\write\gtoutfile{MSC-class: \theprimaryclass\ifx\thesecondaryclass\relax\else, \thesecondaryclass\fi}
\immediate\write\gtoutfile{Journal-ref: Algebr. Geom. Topol. \thevolumenumber\s
(\thevolumeyear) \startpage-\finishpage}
\immediate\write\gtoutfile{Comments: Published by Algebraic and
Geometric Topology at}
\immediate\write\gtoutfile{\s\s\s  http://www.maths.warwick.ac.uk/agt/AGTVol\thevolumenumber/agt-\thevolumenumber-\thepapernumber.abs.html}
\immediate\write\gtoutfile{\noexpand\\}
\immediate\write\gtoutfile{}
\ifx\theasciiabstract\relax
\immediate\write\gtoutfile{\theabstract}\else
\immediate\write\gtoutfile{\theasciiabstract}\fi
\immediate\write\gtoutfile{}
\immediate\write\gtoutfile{\noexpand\\}
\immediate\write\gtoutfile{}
\immediate\closeout\gtoutfile}}  

\def\maketitlepage{\makeagttitle\makeheadfile}
\let\makeshorttitle\maketitlepage
\let\maketitle\maketitlepage

\lognumber{43}
\volumenumber{2}
\volumeyear{2002}
\papernumber{43}
\pagenumbers{1119}{1145}
\received{22 October 2002}
\accepted{30 November 2002}
\published{7 December 2002}

\usepackage{amsmath,amssymb}   
\usepackage{graphicx}
\usepackage{psfrag}
\usepackage{amscd}

\def\psfraga <#1,#2> #3#4{\psfrag #3 
{\smash{\rlap{\kern #1 \raise #2\hbox{#4}}}}}

\newcommand{\Exp}[2]{\ensuremath{\exp_{#1}\!{#2}}}
\newcommand{\s}{\ensuremath{S^1}}
\newcommand{\exps}[1]{\Exp{#1}{\s}}
\newcommand{\expcupone}[1]{\Exp{#1}{(\s,1)}}

\newcommand{\real}{\mathbf{R}}
\newcommand{\integer}{\mathbf{Z}}
\newcommand{\rational}{\mathbf{Q}}
\newcommand{\complex}{\mathbf{C}}

\newcommand{\twobytwo}[4]
{\left[\begin{array}{cc} #1 & #2 \\ #3 & #4 \end{array}\right]}
\newcommand{\intr}{\mathrm{int}\,}

\newtheorem{theorem}{Theorem}
\newtheorem{lemma}{Lemma}
\newtheorem{corollary}{Corollary}
\theoremstyle{remark}
\newtheorem*{remark}{Remark}

\begin{document}

\title{Finite subset spaces of $S^1$}
\asciititle{Finite subset spaces of S^1}
\author{Christopher Tuffley}
\address{Department of Mathematics, University of California\\Berkeley, 
CA 94720, U.S.A.}
\email{tuffley@math.berkeley.edu}


\begin{abstract}
Given a topological space $X$ denote by $\Exp{k}{X}$ the space of 
non-empty subsets of $X$ of size at most $k$, topologised as a quotient 
of $X^k$. This space may be
regarded as a union over $1\leq l \leq k$ of configuration spaces
of $l$ distinct unordered points in $X$.
In the special case
$X=\s$ we show that: (1) \Exp{k}{\s}\ has the homotopy type of an odd
dimensional sphere of dimension $k$ or $k-1$; (2) the natural inclusion of
$\Exp{2k-1}{\s}\simeq S^{2k-1}$ into $\Exp{2k}{\s}\simeq S^{2k-1}$ is
multiplication by two on homology; (3) the complement 
$\Exp{k}{S^1}\setminus\Exp{k-2}{S^1}$ of the codimension two strata in
$\Exp{k}{\s}$ has the homotopy type of a $(k-1,k)$--torus knot complement;
and (4) the degree of an induced map 
$\Exp{k}{f}\co\exps{k}\rightarrow\exps{k}$ is 
$(\deg f)^{\lfloor (k+1)/2\rfloor}$ for $f\co\s\rightarrow\s$.
The first three results generalise known facts that \exps{2}\ is a
M\"obius strip with boundary \exps{1}, and that \exps{3}\ is the three-sphere
with \exps{1}\ inside it forming a trefoil knot.
\end{abstract}

\asciiabstract{Given a topological space X denote by exp_k(X) the space of 
non-empty subsets of X of size at most k, topologised as a quotient of
X^k. This space may be regarded as a union over 0 < l < k+1 of
configuration spaces of l distinct unordered points in X.  In the
special case X=S^1 we show that: (1) exp_k(S^1) has the homotopy type
of an odd dimensional sphere of dimension k or k-1; (2) the natural
inclusion of exp_{2k-1}(S^1) h.e. S^{2k-1} into exp_2k(S^1)
h.e. S^{2k-1} is multiplication by two on homology; (3) the complement
exp_k(S^1)-exp_{k-2}(S^1) of the codimension two strata in exp_k(S^1)
has the homotopy type of a (k-1,k)-torus knot complement; and (4) the
degree of an induced map exp_k(f): exp_k(S^1)-->exp_k(S^1) is (deg
f)^[(k+1)/2] for f: S^1-->S^1.  The first three results generalise
known facts that exp_2(S^1) is a Moebius strip with boundary
exp_1(S^1), and that exp_3(S^1) is the three-sphere with exp_1(S^1)
inside it forming a trefoil knot.}

\maketitle

\section{Introduction}

\subsection{Finite subset spaces}

Given a topological space $X$ let \Exp{k}{X} denote the set of
all nonempty finite subsets of $X$ of cardinality at most $k$. There is a
natural map
\begin{eqnarray*}
\underbrace{X \times \cdots \times X}_k & \rightarrow & \Exp{k}{X} \\
(x_1,\ldots,x_k)            & \mapsto     & \{x_1\}\cup\cdots\cup\{x_k\}
\end{eqnarray*}
and we endow \Exp{k}{X} with the quotient topology to obtain a topological
space, the $k$th finite subset space of $X$. The first finite subset space
\Exp{1}{X}\ is 
clearly $X$ for all $X$, and \Exp{2}{X} co-incides
with the second symmetric product $\mathrm{Sym}^2(X)$, but for $k\geq 3$ we 
have a proper quotient of $\mathrm{Sym}^k(X)$ since, for example, the
points $(a,a,b)$ and $(a,b,b)$ in $X^3$ both map to $\{a,b\}$ in \Exp{3}{X}.
These extra identifications mean that \Exp{k}{X} is in general highly
singular, but give rise to natural inclusions
\begin{eqnarray*}
\Exp{l}{X}  & \hookrightarrow & \Exp{k}{X} \\
\{x_1,\ldots,x_j\} &  \mapsto  & \{x_1,\ldots,x_j\}
\end{eqnarray*}
for $j\leq l\leq k$, maps that require a choice of basepoint for 
$\mathrm{Sym}^k(X)$. The space \Exp{k}{X}
may thus be regarded as a union over $1\leq l \leq k$ of configuration spaces
of $l$ distinct unordered points in $X$. Moreover \Exp{k}{X}\ is compact
whenever $X$ is, in which case it gives a compactification of the corresponding
configuration space. Such spaces and their compactifications 
have been of considerable interest
recently in algebraic topology. See, for example, Fulton and 
MacPherson~\cite{fulton-macpherson94}, Levitt~\cite{levitt95}, 
Yoshida~\cite{yoshida96}, and Ulyanov~\cite{ulyanov02}.

Given a map $f\co X\rightarrow Y$ we obtain a map 
$\Exp{k}{f}\co\Exp{k}{X}\rightarrow\Exp{k}{Y}$ in the obvious way, by sending
$S\subseteq X$ to $f(S)\subseteq Y$. This construction turns $\Exp{k}{}$
into a functor. Moreover, if $\{h_t\}$ is a homotopy between $f$ and $g$
then $\{\Exp{k}{h_t}\}$ is a homotopy between $\Exp{k}{f}$ and $\Exp{k}{g}$,
so that $\Exp{k}{}$ is in fact a functor on the level of homotopy classes of
maps and spaces.

The space \Exp{k}{X}\ was introduced by Borsuk and Ulam~\cite{borsukulam31}
in 1931 as the symmetric product, and has been re-introduced more
recently by Handel~\cite{handel00} in a paper that
establishes many general properties for Hausdorff $X$ and some interesting
homotopy properties when additionally $X$ is path-connected. 
Various different notations have been used for \Exp{k}{X}, including 
$X(k)$, $X^{(k)}$, $\mathcal{F}_k(X)$ and $Sub(X,k)$;
our notation follows that used by 
Mostovoy~\cite{agtmostovoy99} and reflects the idea that we are truncating
the (suitably interpreted) series
\[
\exp X = \emptyset \cup X \cup \frac{X^2}{2!}
\cup \frac{X^3}{3!} \cup \cdots
\]
at the $X^k/k!=X^k/S_k$ term. The name, however, is our own. There does
not seem to be a satisfactory name in use among geometric 
topologists---indeed, recent authors Mostovoy and Handel do not use
any name at all---and while symmetric product has remained in use 
among authors such as Illanes~\cite{illanes85} and Mac\'{\i}as~\cite{macias99}
writing from the perspective of general topology we prefer to use this for 
$X^k/S_k$. We therefore
propose the descriptive name $k$th finite subset space used here. 

In what follows we will be concerned exclusively with the case $X=\s$. 
The results are not only pretty, but also of topological interest 
due to their connection with configuration spaces and their compactifications.

\subsection{Known and new results on \Exp{k}{S^1}}
\label{knownandnew.sec}

A sequence of pictures, outlined in section~\ref{moebius.sec}, shows that 
\exps{2} is a M\"{o}bius strip with boundary \exps{1}. Note in
particular that both \exps{1}\ and \exps{2}\ have the homotopy type of
\s\, and that the inclusion map induces multiplication by two on $H_1$. 
The homeomorphism type of \exps{3} is also known and was calculated by Bott, 
correcting Borsuk's 1949 paper~\cite{borsuk49}:

\begin{theorem}[Bott~\cite{bott52}]
\label{th.Bott}
The space \exps{3} is homeomorphic to the 3-sphere $S^3$.
\end{theorem}

Bott proves this using a cut-and-paste argument, first showing that
\exps{3}\ may be obtained from a single 3-simplex by gluing faces in
pairs, then using this to find $\pi_1(\exps{3})=\{1\}$. He then divides the
simplex into a number of pieces which he re-assembles to form solid tori,
which give $\exps{3}$ when glued along their boundary. This exhibits
$\exps{3}$ as a simply connected lens space, hence $S^3$. An explicit
homeomorphism is not given, and indeed it is non-obvious, as the following 
theorem of E. Shchepin illustrates.

\begin{theorem}[Shchepin, unpublished] 
\label{th.trefoil}
The inclusion $\exps{1} \hookrightarrow \exps{3}$ is a trefoil knot.
\end{theorem}
As further illustration, the two-cells in the above 
simplicial
decomposition of $S^3$ form a M\"obius strip and a dunce cap. The first
is of course \exps{2} bounding \exps{1}, and the second consists of those
subsets containing $1\in \s$.

Shchepin's proof of Theorem~\ref{th.trefoil} is apparently based on 
a direct calculation of the fundamental group of 
$\exps{3}\setminus \exps{1}$. We will give two independent simultaneous proofs
of Theorems~\ref{th.Bott} and~\ref{th.trefoil}, one via cut and paste
topology, and the second via the classification of Seifert fibred spaces.
The natural action of \s\ on itself gives an action of \s\ on \exps{k}\
for each $k$, and for $k=3$ we obtain the following refinement of
Bott's and Shchepin's results:

\begin{theorem}
\label{th.seifertfibred}
The space \exps{3}\ is a Seifert fibred $3$--manifold, and as such is
oriented fibre preserving diffeomorphic to $S^3$ with the $(2,-3)$ \s\ action
\begin{equation}
\lambda \cdot (z_1,z_2) = (\lambda^2 z_1, \lambda^{-3} z_2) ,
\label{23action.eq}
\end{equation}
where we regard \s and $S^3$ as sitting in $\complex$ and 
$\complex^2$ respectively and give \exps{3}\ the canonical orientation
it inherits from \s.
\end{theorem}

We mention also an elegant geometric construction due to 
Mostovoy~\cite{agtmostovoy99} showing that both 
Theorems~\ref{th.Bott} and~\ref{th.trefoil} can be deduced from
known facts about lattices in the plane. Namely, the space 
$SL(2,\real)/SL(2,\integer)$ of plane lattices modulo scaling
is diffeomorphic to a trefoil complement (the proof, due to D. Quillen,
may be found in~\cite[page 84]{milnorsktheory} and uses the Weierstrass $\wp$ 
function associated to the lattice), and this space may be compactified
by adding degenerate lattices to obtain $S^3$. 
Together with Theorems~\ref{th.Bott} and~\ref{th.trefoil} this shows
$\exps{3}\setminus \exps{1} \cong SL(2,\real)/SL(2,\integer)$, and
Mostovoy's construction fills in the the third side of this triangle, 
associating to each lattice a finite subset of \s, thought of as 
$\real P^1$. Each degenerate lattice corresponds to a one element
subset, each rectangular lattice a two element subset, and all other lattices
correspond to three element subsets. Moreover, his map is equivariant with
respect to the natural actions of $\s$ on \exps{3} and 
$PSO(2) \subseteq PSL(2,\real)$ on $SL(2,\real)/SL(2,\integer)$.

The first three of the following new results generalise the theorems and 
observations above. Proofs appear in subsequent sections. 
Since writing this paper I have learnt that
Theorems~\ref{th.odddimsphere},~\ref{th.timestwo} and the observation
on the map $\exps{k-1}\rightarrow\exps{k}:\Lambda\mapsto\Lambda\cup\{1\}$
that follows their proofs have
been proved independently by 
David Handel~\cite[unpublished work]{handelcircle}, using essentially the
same argument.

\begin{theorem}
\label{th.odddimsphere} The space
\exps{k}\ has the homotopy type of an odd dimensional sphere, of dimension
$k$ or $k-1$ according to whether $k$ is odd or even.
\end{theorem}

Since $\exps{2k-1}\simeq S^{2k-1} \simeq \exps{2k}$ we may ask how 
\exps{2k-1}\ sits inside \exps{2k}. The following result falls out of
the proof of Theorem~\ref{th.odddimsphere} and shows that the situation
is analogous to \exps{1}\ inside \exps{2}:

\begin{theorem}
The inclusion $\exps{2k-1} \hookrightarrow \exps{2k}$ induces 
multiplication by two on $H_{2k-1}$.
\label{th.timestwo}
\end{theorem}

As our last generalisation, \exps{k-2}\ inside \exps{k} is in some sense a 
homotopy $S^{l-2}$ inside a homotopy $S^l$, so it is natural to ask if this is 
embedded in some interesting way too. Analogously to Theorem~\ref{th.trefoil} 
we have:

\begin{theorem}
\label{th.torusknot}
The complement of \exps{k-2}\ in \exps{k}\ has the homotopy type of
a $(k-1,k)$--torus knot complement.
\end{theorem}

Finally, a map $f\co\s \rightarrow \s$ induces a map $\Exp{k}{f}$ of 
homotopy spheres, and we calculate its degree in terms of $k$ and the degree
of $f$.

\begin{theorem}
\label{th.degree}
If $f\co\s\rightarrow\s$ then 
\[
\deg \Exp{k}{f} = \left(\deg f\right)^{\left\lfloor\frac{k+1}{2}\right\rfloor}.
\]
\end{theorem}

\begin{remark}
It is perhaps worth noting that although \exps{k}\ is a manifold for
$k=1,2,3$, it is not a manifold for any $k\geq 4$. 
Subsets of size $k$ do have $k$--ball neighbourhoods in \exps{k}\ and this
transition may be understood in terms of neighbourhoods of $k-1$ element
subsets as follows. Given such a subset $\Lambda$ of \s\ there are $k-1$
choices for which point to consider ``doubled'' and split in two to
obtain a $k$ element subset. Each such choice leads to a $k$--dimensional
halfspace containing $\Lambda$ and we obtain a neighbourhood of $\Lambda$ by 
gluing these together along their boundaries. The transition thus occurs when 
$k-1$ increases above two, and we see also that when $k=2$ points in \exps{1}\ 
have halfspace neighbourhoods and thus form the boundary of a $2$--manifold.

Later we shall see this more explicitly when we show that \exps{k}\ may be 
obtained from a single $k$--simplex by identifying faces. Under the
identifications the $0$th and $k$th faces become one face and the remaining
$k-1$ faces a second, corresponding to \exps{k-1}.
\end{remark}

\subsection{Notation and terminology}

The proofs of our main results make use of arguments involving simplices
and simplicial decompositions and we take a moment to fix language. For
the most part we follow Hatcher~\cite[section 2.1]{hatcher}, and will use
``simplicial decomposition'' in the sense of his 
$\Delta$--complexes.
In particular, we will not require that the simplices in our decompositions
be determined by their vertices.

Given $k+1$ points $u_0,\ldots,u_k$ in an affine space let
\[
[u_0,\ldots,u_k] = \left\{ \mbox{$\sum_i t_i u_i$} \big| 
   \mbox{$\sum_i t_i = 1$ and $t_i\geq 0$ for all $i$} \right\},
\]
the set of convex combinations of the $u_i$. We will typically only
write this when the points are affinely independent, in which case
$[u_0,\ldots,u_k]$ is a $k$--simplex with (ordered) vertices $u_0,\ldots,u_k$.
For $\sigma\in S_{k+1}$ we regard $[u_{\sigma(0)},\ldots,u_{\sigma(k)}]$ as 
the same simplex with orientation $(-1)^{\mathrm{sign}\,\sigma}$ times that of 
$[u_0,\ldots,u_k]$.

The canonical map between two simplices $[u_0,\ldots,u_k]$ and 
$[v_0,\ldots,v_k]$ is the unique map given by sending $u_i$ to $v_i$ for 
each $i$ and extending affinely. A hat~$\hat{~}$ over a vertex means it 
is to be omitted, in other words
\[
[u_0,\ldots,\hat{u}_{i},\ldots,u_k]= [u_0,\ldots,u_{i-1},u_{i+1},\ldots,u_k],
\]
and unless indicated otherwise the interior of a simplex $[u_0,\ldots,u_k]$
will always mean the open simplex
\[
\intr [u_0,\ldots,u_k] = \left\{ \mbox{$\sum_i t_i u_i$} \big| 
   \mbox{$\sum_i t_i = 1$ and $t_i> 0$ for all $i$} \right\}
\]
regardless of whether this is an open subset of the ambient space.

\section{Finite subset spaces of \s\ have the homotopy type of 
odd dimensional spheres}
\label{sec.odddimsphere}

\subsection{Introduction}

To prove that \exps{k}\ has the homotopy type of a sphere we will find a
cell structure for it and use this to show it has the correct fundamental
group and homology. Application of a standard argument that a simply
connected homology sphere is a homotopy sphere then yields the result. Before
doing so however, let us look at \exps{2}\ and \exps{3}\ in some detail,
which will illustrate the situation in higher dimensions.

\begin{figure}[b]
\begin{center}\small
\leavevmode
\psfrag{(x,y)}{$(x,y)$}
\psfrag{(y,x)}{$(y,x)$}
\psfrag{(a)}{(a)}
\psfrag{(b)}{(b)}
\psfrag{(c)}{(c)}
\psfrag{(d)}{(d)}
\includegraphics{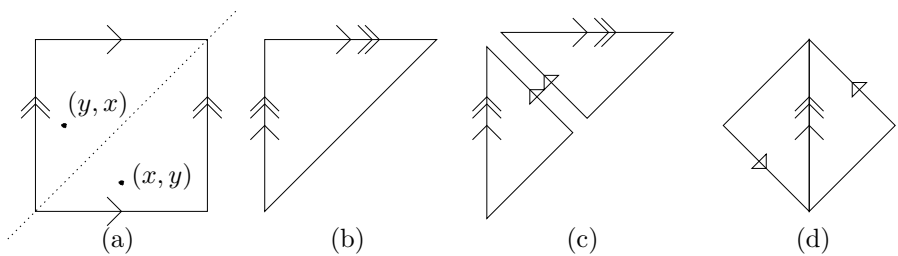}
\caption{The space \exps{2}\ is a M\"obius strip. 
         Fold along the diagonal in (a)
         to identify $(x,y)$ and $(y,x)$; the result is shown in (b). 
  	 Cut (c) and re-glue (d) to recognise this as a M\"obius strip.}
\label{exp2s1.fig}
\end{center}
\end{figure}

\subsection{The homeomorphism type of \exps{2}}
\label{moebius.sec}

To see that \exps{2}\ is a M\"obius strip we will use the usual picture
of $\s \times \s$ as a square with opposite sides identified. In forming
\exps{2}\ the point $(x,y)$ is identified with $(y,x)$, so the square is
folded along the diagonal shown dotted in figure~\ref{exp2s1.fig}(a),
resulting in the triangle with two edges identified shown in
figure~\ref{exp2s1.fig}(b). This possible unfamiliar picture may be
recognised as a M\"obius strip either by cutting and re-gluing as in
figures~\ref{exp2s1.fig}(c) and~(d), or by gluing just the ends of the
hypotenuse together to get a punctured projective plane.

\begin{figure}[ht!]
\begin{center}
\leavevmode
\includegraphics{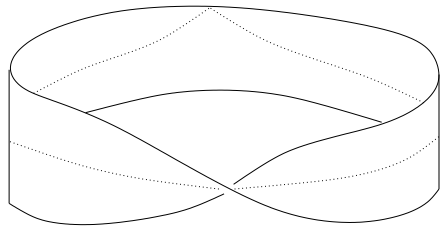}
\caption{A M\"obius strip, in which the glued edge-pair of the triangle of
         figure~\ref{exp2s1.fig}(b) corresponding to the set 
        $\{\Lambda \in \exps{2} | 1\in\Lambda\}$ is shown dotted.}
\label{moebius.fig}
\end{center}
\end{figure}

A M\"obius strip with the edge corresponding to the glued sides of the
triangle shown dotted appears in figure~\ref{moebius.fig}. The diagonal
edge forms the boundary circle. Note that the diagonal comes from
the diagonal embedding of $S^1$ in $\s \times \s$, which maps to 
$\exps{1}\subseteq\exps{2}$. 

We mention also a construction pointed out by Chuck Livingston. There is
a natural map from \exps{2}\ to $\real P^1$ sending each pair of points
on \s\ to the line through the origin bisecting the arc between them. 
The fibre above each point in $\real P^1$ is an interval of length $\pi$ 
and the bundle is easily seen to be non-orientable.

\subsection{The homeomorphism types of \exps{3}\ and 
$\exps{3}\setminus\exps{1}$}
\label{exp3.sec}

\begin{figure}[ht!]
\begin{center}\small
\leavevmode
\psfrag{v0}{$v_0$}
\psfrag{v1}{$v_1$}
\psfrag{v2}{$v_2$}
\psfrag{v3}{$v_3$}
\psfrag{x}{$x$}
\psfrag{y}{$y$}
\psfrag{z}{$z$}
\psfrag{(a)}{(a)}
\psfrag{(b)}{(b)}
\psfrag{a}{$a$}
\psfrag{b}{$b$}
\includegraphics{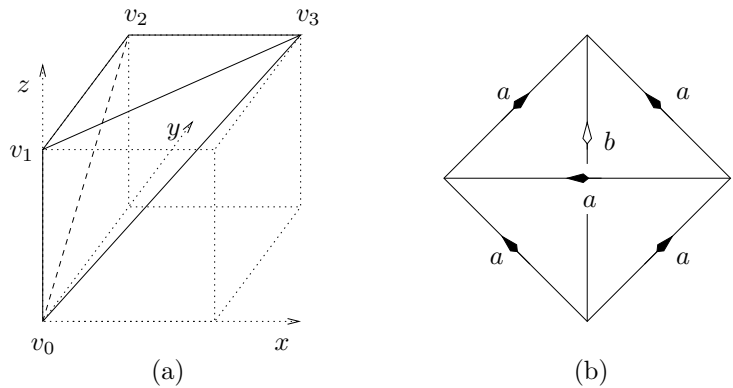}
\caption{Simplicial decomposition of \exps{3}. \exps{3} may be formed 
 	from the simplex $[v_0,v_1,v_2,v_3]$ in (a) by identifying the face
	$[v_0,v_1,v_2]$ with $[v_1,v_2,v_3]$, and the face
	$[v_0,v_1,v_3]$ with $[v_0,v_2,v_3]$. This corresponds 
	to the edge and associated face gluings of (b).}
\label{3simplex.fig}
\end{center}
\end{figure}

Now consider \exps{3}. Begin again with $[0,1]^3=I^3$ with opposite faces
identified. Each $\Lambda\in\exps{3}$ has at least one representative 
$(x,y,z)\in I^3$ with $0\leq x \leq y \leq z \leq 1$, so we may restrict
our attention to the simplex with vertices $v_0=(0,0,0)$, $v_1=(0,0,1)$, 
$v_2=(0,1,1)$ and $v_3=(1,1,1)$ shown in figure~\ref{3simplex.fig}(a). Now
$(0,x,y)\sim (x,y,1)$ in \exps{3}, so the face $[v_0,v_1,v_2]$ is
glued to the face $[v_1,v_2,v_3]$; next $(x,x,y)\sim (x,y,y)$ in \exps{3},
so the face $[v_0,v_1,v_3]$ is glued to $[v_0,v_2,v_3]$. This accounts
for all the identifications of the simplex arising from \exps{3}, and the 
result, taking account of edge identifications, is shown in 
figure~\ref{3simplex.fig}(b), viewed from a different angle. There is just
one vertex (the set $\{1\}$), two edges and 
two $2$--simplices, one forming a dunce cap and the second a M\"obius
strip. The dunce cap comes from the two faces $x=0$ and $z=1$ and corresponds
to the set $\{\Lambda \in \exps{3} | 1\in\Lambda\}$; the M\"obius strip comes
from the faces $x=y$ and $y=z$ and is
of course $\exps{2}\subseteq\exps{3}$.

From above we have $\chi(\exps{3})=1-2+2-1=0$; a well known but still magical
fact (see~\cite[p. 122]{thurston97}) then 
guarantees that the space 
formed by gluing the faces of the 
$3$--simplex in this fashion is a $3$--manifold. Calculating $\pi_1(\exps{3})$
using figure~\ref{3simplex.fig}(b) we obtain the presentation 
$\langle a,b | a, a^2 b^{-1} \rangle$, so \exps{3}\ is in fact a simply
connected $3$--manifold---and hence almost certainly $S^3$, especially given
its simple construction. Bott~\cite{bott52} completes the proof by showing
it has a genus one Heegaard splitting and appealing to the classification
of lens spaces. A more informative approach is the proof of
Theorem~\ref{th.seifertfibred} given below. In 
appendix~\ref{cutandpaste.apdx} (beginning page~\pageref{cutandpaste.apdx}) 
we will also show directly that 
\exps{3}\ is a $3$--sphere (and that \exps{1}\ inside it is a trefoil)
using the genus two Heegaard splitting obtained from the triangulation
above.

Before proving Theorem~\ref{th.seifertfibred} a word is necessary on 
orientation. Given
an orientation of \s\ and a set $\Lambda=\{\lambda,\mu,\nu\}$ of three 
distinct points 
in \s\ we may canonically orient 
\[
T_\Lambda\exps{3}\cong T_\lambda\s \oplus T_\mu\s \oplus T_\nu\s
\]
by positively orienting each of the summands and requiring that the
cyclic ordering of $\lambda$, $\mu$, and $\nu$ agree with the orientation 
of \s . This extends to an orientation of \exps{3}\ and we
regard this as the canonical orientation of \exps{3}. Note that
the standard orientation of the simplex $0\leq x\leq y\leq z\leq 1$ of 
figure~\ref{3simplex.fig} co-incides with the canonical orientation
of \exps{3}.

To orient $S^3$ we regard it as the boundary of the $4$--ball in 
$\complex^2$ with its canonical orientation and use the ``outward first''
convention for induced orientations on boundaries.

\begin{proof}[Proof of Theorem~\ref{th.seifertfibred}]
We have seen that \exps{3}\ is a closed simply connected $3$--manifold and we 
observe that it is Seifert fibred by the natural action of \s. There are 
precisely two exceptional fibres, the orbits of $\{1,-1\}$ and 
$\{1,e^{\pi i/3},e^{2\pi i/3}\}$, and these have multiplicities $2$ and $3$ 
respectively. Since $S^3$ with the $(2,-3)$ action of 
equation~(\ref{23action.eq}) on 
page~\pageref{23action.eq} shares these properties our aim will be to 
show that they are enough to completely determine the unoriented fibre type
of \exps{3}.

To this end let $M$ be a closed simply connected Seifert fibred space with
precisely two exceptional fibres, of multiplicities two and three. Simple
connectivity of $M$ implies the orbit surface is simply connected also, and
therefore $S^2$. Moreover, the fibres may be consistently oriented. Removing
fibred solid torus neighbourhoods of each of the exceptional fibres
thus leaves an oriented circle bundle over a twice punctured sphere, which 
we may write as $\s\times I\times\s$, the base corresponding to the first
two factors. 

$M$ is completely determined by specifying slopes 
$\alpha_0/\beta_0,\alpha_1/\beta_1\in\rational$ along which to glue back in
meridional discs to $\s\times\{0\}\times\s$ and $\s\times\{1\}\times\s$.
To get the correct multiplicities we must have 
$\{\beta_0,\beta_1\}=\{2,3\}$ so without loss of generality let
$\beta_0=2$, $\beta_1=3$. The classification of orientable Seifert
fibred spaces (see for example Hatcher's $3$--manifold 
notes~\cite[p. 25]{hatcher3M}) tells us that
the slopes are only determined $\bmod$ $1$ subject to their sum being
fixed (this ambiguity comes from the choice of trivialisation of the
circle bundle) so we may further assume $\alpha_0=1$ and write 
$\alpha_1=\alpha$. We calculate $\pi_1$ of the resulting manifold. 

The fundamental group of $\s\times I\times\s$ is free abelian, generated by
$b$ and $f$, where $b$ and $f$ are positively oriented generators of
$\pi_1$ of the base and fibre respectively. Gluing a disc in to 
$\s\times\{0\}\times\s$ along a line of slope $1/2$ kills
$2b+f$ while gluing a disc in to $\s\times\{1\}\times\s$ along
a line of slope $\alpha/3$ kills $-3b+\alpha f$, the minus sign coming from the
fact that $\s\times\{1\}\times\s$ has orientation $[-b,f]$ (recall that
we are using ``outward first'' to orient boundaries). Thus simple
connectivity of $M$ implies
\[
\det \twobytwo{2}{-3}{1}{\alpha} = 2\alpha +3 = \pm 1.
\]
Setting the determinant equal to $1$ and $-1$ in turn gives $\alpha=-1$, 
$\alpha=-2$. Thus there are exactly two 
possibilities for the oriented fibre type of $M$; in Hatcher's notation
they are $M(\pm g,b;\alpha_0/\beta_0,\alpha_1/\beta_1)=M(0,0;1/2,-1/3)$
and $M(0,0;1/2,-2/3)$, where $\pm g$ specifies the genus and orientability
of the orbit surface and $b$ the number of boundary components. 

Reversing the orientation of 
$M(\pm g,b;\alpha_1/\beta_1,\ldots,\alpha_k/\beta_k)$ simply changes the
signs of all the attaching slopes. Thus
\[
-M(0,0;1/2,-1/3)\cong M(0,0;-1/2,1/3) \cong M(0,0;1/2,-2/3)
\]
and the unoriented fibre type of $M$ is completely determined as claimed.
It follows that $\exps{3}$ and $S^3$ with the $(2,-3)$ action are 
fibre preserving diffeomorphic and all that remains is to determine whether
orientation is preserved or reversed.

To determine orientations we look at the return map on a disc $D$ transverse
to the exceptional fibre of multiplicity three at a point $p$.
The fibre $F$ is oriented by the \s\ action and we orient $D$ such that
$T_p F\oplus T_p D$ is positive.  In the case of \exps{3}\ we use the
point $p=(1/3,2/3,1)\in I^3$ and a small disc containing it in the plane
$z=1$. The vector $(1,1,1)$ forms a positive basis for $T_p F$ so 
$T_p D$ has orientation $[e_x,e_y]$, and the
action of the return map is the $1/3$ anti-clockwise twist given by
the canonical map from $[v_1,v_2,v_3]$ to $[v_3,v_1,v_2]$.
For $S^3$ we take $p=(0,1)$ and consider the $(2,3)$ action. The tangent space
$T_p S^3$ has orientation $[(0,i),(1,0),(i,0)]$, in which $(0,i)$ forms
a positive basis for $T_p F$, so for a suitable choice of disc $D$ its
tangent space at $p$ has
positive basis $\{(1,0),(i,0)\}$. The first return is when 
$\lambda=\omega:=e^{2\pi i/3}$ and we see that the derivative of the first
return map is multiplication by $\omega^2$ on $T_p D$. This is a clockwise
rotation through $2\pi/3$, so to match orientations with \exps{3}\ we
must reverse the orientation of the orbit through $(0,1)$, giving the
$(2,-3)$ \s\  action as claimed.
\end{proof}

\subsection{The homotopy type of \exps{k}}
\label{expkS1.sec}

Finally we turn our attention to the general case. Proceeding 
analogously to the two and three dimensional cases we may reduce to the
$k$--simplex $0\leq x_1 \leq \ldots \leq x_k \leq 1$. Working somewhat
more formally than above let
\[
v_i = (\,\underbrace{0,\ldots,0}_{k-i},\underbrace{1,\ldots,1}_i\,) \in \real^k
\]
for $i=0,\ldots,k$, and let $\sigma_k$ be the map of the simplex
$[v_0,\ldots,v_k]$ to \exps{k}, $\tau_{k-1}$ that of the simplex 
$[v_1,\ldots,v_k]$. Being a little sloppy with notation we claim:

\begin{lemma}
\label{lem.cellstructure}
\exps{k} has a simplicial decomposition with one 0-simplex, two $i$--simplices
for each $1\leq i \leq k-1$, and one $k$--simplex, namely, $\sigma_i$ for
$0\leq i \leq k$ and $\tau_i$ for $1\leq i \leq k-1$. The boundary
map $\partial_i \co \integer\sigma_i \oplus \integer\tau_i \rightarrow
\integer\sigma_{i-1} \oplus \integer\tau_{i-1}$ has matrix
\[
D_i = \frac{1+(-1)^i}{2}\twobytwo{-1}{0}{2}{1}
\]
for $2\leq i\leq k-1$, is the zero map for $i=1$ and has matrix
$D_k |_{\integer\sigma_k}$
for $\partial_k \co \integer\sigma_k\rightarrow
\integer\sigma_{k-1} \oplus \integer\tau_{k-1}$, $k\geq 2$.
\end{lemma}

The lemma enables us to calculate the homology of \exps{k}, obtaining
the following.

\begin{corollary}
\label{cor.homology} The space
$\exps{k}$ has the homology of an odd dimensional sphere, of dimension
$k$ if $k$ is odd and dimension $k-1$ if $k$ is even. The inclusion map
$\exps{2k-1}\hookrightarrow\exps{2k}$ induces multiplication by two on
$H_{2k-1}$.
\end{corollary}

\begin{proof}[Proof of Lemma~\ref{lem.cellstructure} and 
Corollary~\ref{cor.homology}]
The lemma is proved by induction. Taking $k=1$ as our base case (although the
cases $k=2$ and~$3$ are largely established above) this is just the usual
cell decomposition of $S^1$ as $I/\{0,1\}$, so consider $k\geq 2$. $\sigma_k$
maps $[v_0,\ldots,v_k]$ onto \exps{k}, taking the interior of the 
simplex homeomorphically 
onto its image, so we need only sort out the face gluings. Faces of the 
form $[v_0,\ldots,\hat{v}_{i},\ldots,v_k]$ where $1\leq i \leq k-1$ 
correspond to 
$0\leq x_1 \leq \ldots \leq x_{k-i}=x_{k-i+1} \leq \ldots \leq x_k\leq 1$, 
giving
subsets of $S^1$ of size $k-1$ or less; therefore $\sigma_k$ restricted to
such a simplex factors through $\sigma_{k-1}$. More precisely, we should
note that the map of simplices factoring this restriction preserves 
orientation. The simplex $[v_0,\ldots,v_{k-1}]$ is
$0=x_1\leq\ldots\leq x_k\leq 1$ which is identified with 
$0\leq x_1\leq\ldots\leq x_k=1$, so $\sigma_k |_{[v_0,\ldots,v_{k-1}]}$
factors through $\tau_{k-1}$ via the canonical map 
$[v_0,\ldots,v_{k-1}]\rightarrow [v_1,\ldots,v_k]$. Thus
\begin{eqnarray*}
\partial \sigma_k & = & \sum_{i=0}^k 
                    (-1)^i  \sigma_k |_{[v_0,\ldots,\hat{v}_{i},\ldots,v_k]} \\
 & = & \tau_{k-1} + \sum_{i=1}^{k-1} (-1)^i \sigma_{k-1} + (-1)^k \tau_{k-1} \\
 & = & \frac{1+(-1)^k}{2} ( -\sigma_{k-1} + 2\tau_{k-1} ),
\end{eqnarray*}
giving the first column of $D_k$.

Turning our attention now to $\tau_{k-1}$, observe that it maps 
$[v_1,\ldots,v_k]$ onto $\{\Lambda \in \exps{k} | 1\in\Lambda \}$, taking
the interior of this simplex homeomorphically onto its image. A
face of $[v_1,\ldots,v_k]$ corresponds to replacing an inequality with 
an equality in
$0\leq x_1 \leq \ldots \leq x_{k-1} \leq x_k = 1$,
and in each case 
maps onto the $k-1$ or fewer element subsets containing $1$. Thus $\tau_{k-1}$
factors through $\tau_{k-2}$ (via an orientation preserving map of simplices)
when restricted to each face, and
\begin{eqnarray}
\partial \tau_{k-1} & = & \sum_{i=1}^{k} 
    (-1)^{i-1}  \tau_{k-1} |_{[v_1,\ldots,\hat{v}_{i},\ldots,v_k]} \nonumber \\
 & = & \sum_{i=1}^{k} (-1)^{i-1} \tau_{k-2} \nonumber \\
 & = & \frac{1+(-1)^{k-1}}{2}\tau_{k-2} , \label{eq.boundarytau} 
\end{eqnarray}
giving the second column of $D_{k-1}$, the first coming from the
inductive hypothesis. This establishes the lemma.

It is now a simple matter to calculate the homology of \exps{k}. The matrix
$D_i$ is zero if $i$ is odd and has determinant $-1$ if $i$ is even, so the
chain maps are alternately zeroes and isomorphisms in the middle dimensions.
Thus $H_i(\exps{k}) = \{0\}$ for $1\leq i \leq k-2$. Clearly 
$H_0(\exps{k})\cong\integer$, so it remains to determine only $H_{k-1}$ and
$H_k$. When $k$ is odd the top end of the chain complex is
\[
0 \longrightarrow \integer 
  \stackrel{0}{\longrightarrow}\integer \oplus \integer
  \stackrel{\cong}{\longrightarrow}\integer\oplus\integer
  \longrightarrow 0 ,
\]
so $H_{k-1}$ is zero and $H_k(\exps{k})\cong\integer$, generated
by $[\sigma_k]$. When $k$ is even we have instead
\[
\begin{array}{c}
0  \longrightarrow  \integer  \longrightarrow  \integer\oplus\integer 
\longrightarrow  0 , \\
                       1       \longmapsto \!     (-1,2)            
\end{array}
\]
making $H_k$ zero and imposing the relation $\sigma_{k-1}=2\tau_{k-1}$
on $\ker \partial_{k-1} = C_{k-1}$. Thus for $k$ even we have 
$H_{k-1}(\exps{k})\cong\integer$, generated by $[\tau_{k-1}]$. This proves
the first statement of the corollary, and for the second simply observe
that the generator $[\sigma_{2k-1}]$ of $H_{2k-1}(\exps{2k-1})$ is twice
$[\tau_{2k-1}]$, the generator of $H_{2k-1}(\exps{2k})$.
\end{proof}

Theorem~\ref{th.odddimsphere} for $k\geq 3$ now follows from 
Corollary~\ref{cor.homology} and the simple connectivity of \exps{3}\
by an application of a standard argument. 

\begin{corollary}[Theorem~\ref{th.odddimsphere}]
\label{cor.odddimsphere}
\exps{k}\ has the homotopy type of an odd dimensional sphere of dimension
$d_k=2\lceil k/2 \rceil -1$.
\end{corollary}

\begin{proof}
We have already seen this for $k=1,2$ so we may assume $k\geq 3$. Then the
$2$--skeleton of \exps{k}\ co-incides with the $2$--skeleton of \exps{3}, 
which is simply connected, so $\pi_1 (\exps{k})\cong \{1\}$ too. By
the Hurewicz theorem $\pi_{d_k}(\exps{k})\cong H_{d_k} (\exps{k})\cong
\integer$, so let $\phi \co S^{d_k} \rightarrow \exps{k}$ be a generator
for $\pi_{d_k}(\exps{k})$. $\phi_*$ induces an isomorphism on $\pi_{d_k}$
and so on $H_{d_k}$ also; since $H_0$ and $H_{d_k}$ are the only 
non-vanishing homology groups of both $S^{d_k}$ and \exps{k}\, $\phi_*$
is an isomorphism on $H_n$ for all $n$. By the simply connected version
of Whitehead's theorem that only requires isomorphisms on homology, $\phi$
is a homotopy equivalence.
\end{proof}

We close this section with some remarks on two related spaces. Let
\begin{eqnarray*}
\expcupone{k} & = & \bigl\{\Lambda \in \exps{k} \big| 1\in \Lambda \bigr\} \\
              & = & \bigl\{\Lambda\cup\{1\}\big| \Lambda \in \exps{k-1}\bigr\},
\end{eqnarray*}
the subsets of \s\ of size $k$ or less that contain $1$. This subspace
has a cell structure with one cell $\tau_i$ in each dimension less than or 
equal to $k-1$, and by (\ref{eq.boundarytau}) the boundary maps are alternately
zero and isomorphisms. Thus the reduced homology of \expcupone{k}\ is 
given by
\[
\hat{H}_i(\expcupone{k}) = \left\{
\begin{array}{cl}
\integer & \mbox{if $i=k-1$ is odd} \\
0        & \mbox{otherwise,}
\end{array}\right. \;
\]
and moreover \expcupone{k}\ is simply connected for $k\geq 3$ since its
2-skeleton is a dunce cap. It follows (by the Whitehead theorem for $k$
odd, and the argument of Corollary~\ref{cor.odddimsphere} for $k$ even)
that \expcupone{k}\ is contractible if $k$ is odd, and homotopy equivalent
to $S^{k-1}$ if $k$ is even. The natural map 
\[
\exps{k-1}\rightarrow\expcupone{k} : \Lambda\mapsto\Lambda\cup\{1\}
\]
takes $[\sigma_{k-1}]$ to $[\tau_{k-1}]$ and so is a homotopy equivalence for
$k$ even. Note however that this map is not a homeomorphism except when
$k=2$ since otherwise $\{\lambda\}$ and $\{\lambda,1\}$ have the same image.

Lastly let
\[
\Exp{}{\,\s} = \bigcup_{k=1}^\infty \exps{k}=  
\{ \Lambda \subseteq \s \big| 0< |\Lambda|<\infty \},
\]
topologised as the direct limit of the \exps{k}, or equivalently, with the
CW topology coming from the cell structure consisting of the 
$\sigma_i,\tau_i$. The full finite subset space $\Exp{}{\,\s}$ has vanishing 
reduced homology in all dimensions and is therefore contractible.

\section{Removing the codimension two strata gives a homotopy torus knot 
complement}

\subsection{A new model for \exps{k}}
\label{newmodel.sec}

\begin{figure}[t]
\begin{center}\small
\leavevmode
\psfrag{a}{$a$}
\psfrag{b}{$b$}
\psfrag{c}{$c$}
\psfrag{ak}{$a^{k-1}$}
\psfrag{bk}{$b^k$}
\includegraphics{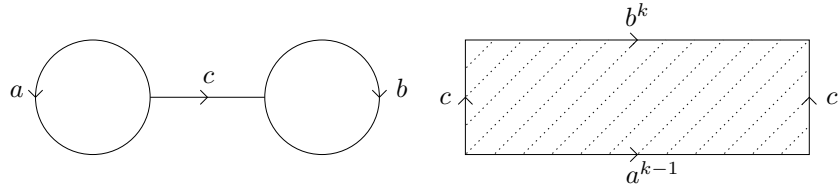}
\caption{A cell structure consisting of three $1$--cells $a,b,c$ and
         one $2$--cell, attached along $a^{k-1} c b^{-k} c^{-1}$.}
\label{cellstructure.fig}
\end{center}
\end{figure}

\begin{figure}[b]
\begin{center}
\leavevmode
\includegraphics{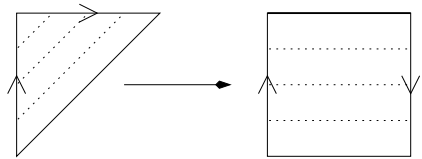}
\caption{The new model in the case $k=2$. 
We model \exps{2}\ as a quotient of a $1$--simplex cross an
interval by normalising the dotted lines in the triangle on the left
(corresponding to subintervals of the \s\ action) to have length one, 
giving the square on the right. In doing this the top left vertex of the 
triangle is stretched to become the top edge of the square.}
\label{simplexCrossI.fig}
\end{center}
\end{figure}

The strategy for proving Theorem~\ref{th.torusknot} is to show that
$\exps{k}\setminus\exps{k-2}$ deformation retracts to a subspace having
the cell structure shown in figure~\ref{cellstructure.fig}. In order to 
make calculations easier and the result more transparent we will
adopt a slightly different picture of \exps{k}, using the action of
\s\ to model it as a quotient of a $(k-1)$--simplex cross an interval.
On the simplex level this action corresponds to the
constant vector field equal to $(1,1,\ldots,1)$ everywhere and the
orbits are unions of intervals of the form
\begin{equation}
\left[(0,x_2,\ldots,x_k)\, ,\,(y_1,\ldots,y_{k-1},1)\right],
\label{eq.S1intervals}
\end{equation}
where $y_i=x_i+1-x_k$  for each $i$. 
The effect of our new model will be to 
normalise the lengths of these intervals to one. 
Figure~\ref{simplexCrossI.fig} depicts this in the case $k=2$.

Let $\Delta$ be the simplex 
\[
[v_0,\ldots,v_{k-1}] = 
\{(0,a_1,\ldots,a_{k-1})|0\leq a_1 \leq \ldots \leq a_{k-1}\leq 1\}
\]
and consider the map from $\Delta\times I$ to our usual model 
\[
X=[v_0,\ldots,v_k]=
\{(x_1,\ldots,x_k)|0\leq x_1 \leq \ldots \leq x_k\leq 1\}
\]
given by
\[
x_i = \left\{\begin{array}{rl}
(1-a_{k-1})t  & i=1, \\
 a_{i-1} + (1-a_{k-1})t & 2\leq i \leq k.
\end{array}\right.
\]
This has inverse 
\begin{eqnarray*}
a_i & = & x_{i+1}-x_1 \\
t   & = & \frac{x_1}{1-x_k+x_1}
\end{eqnarray*}
well defined off the codimension two face $[v_1,\ldots,v_{k-1}]$, which has 
preimage the entire codimension one face $[v_1,\ldots,v_{k-1}]\times I$
due to the intervals~(\ref{eq.S1intervals}) through these points being
stretched from length zero to one. We will refer to this as 
the ``fake face'' of $\Delta\times I$ and denote the quotient map
$\Delta\times I\rightarrow X \rightarrow \exps{k}$ by $q$.

In forming \exps{k}\ from $\Delta\times I$ the $k-1$ faces $\{a_1=0\}\times I$ 
and $\{a_i=a_{i+1}\}\times I$, $1\leq i\leq k-2$, are all identified
according to the maps
\begin{equation}
\begin{array}{ccc}
[v_0,\ldots,\hat{v}_{i},\ldots,v_{k-1}]\times I \rightarrow
      [v_0,\ldots,\hat{v}_{j},\ldots,v_{k-1}]\times I & & (i,j\not=0)
\end{array}
\label{facegluings.eq}
\end{equation}
given by the product of the canonical map with the identity. The
face $\{a_{k-1}=1\}\times I$ is collapsed back down to $\{a_{k-1}=1\}$ by
projection on the first factor, and $\Delta\times\{0\}$ is glued to 
$\Delta\times\{1\}$ according to $(a,1)\sim(\phi(a),0)$, where 
\[
\phi_i(a) = \left\{\begin{array}{ll}
1-a_{k-1}  & i=1, \\
 a_{i-1} + 1-a_{k-1} & 2\leq i \leq k.
\end{array}\right.
\]
$\phi$ is affine and permutes the vertices $v_0,\ldots,v_{k-1}$ 
cyclicly according
to the permutation $i\mapsto i-1 \pmod{k-1}$ and so is the canonical map
$[v_0,\ldots,v_{k-1}]\rightarrow [v_{k-1},v_0,\ldots,v_{k-2}]$. In particular
$\exps{k}\setminus\exps{k-1}$, as the quotient of 
$(\intr\Delta)\times I$, is the mapping torus of 
$\phi|_{\intr \Delta}$ and has the homeomorphism
type of an $\real^{k-1}$ bundle over \s, with monodromy of order $k$.
$\phi$ reverses orientation exactly when $k$ is even so the bundle
is trivial for $k$ odd and nontrivial for $k$ even.

Although we shall not explicitly do so, there is no loss in generality in 
regarding $\Delta$ as the more symmetrical standard $(k-1)$--simplex
\[
\bigl\{ (t_0,t_1,\ldots,t_{k-1})\in\real^k \big| 
          \mbox{$\sum_i t_i = 1$ and $t_i\geq 0$ for all $i$} \bigr\},
\]
and what follows may be read with this picture in mind.

\subsection{The fundamental group of $\exps{k}\setminus\exps{k-2}$}

As a first application of our new model we calculate the fundamental group
of $\exps{k}\setminus\exps{k-2}$. This calculation is in some sense
redundant, in that in the following section we will find a $2$--complex
to which it is homotopy equivalent. However, the proof given below
that this space has the correct fundamental group strongly echoes
the corresponding calculation for the torus knot complement. In
doing so it carries the main insight as to why the two have the same
homotopy type, while the proof of this fact, while geometric
in nature, is somewhat technical. We therefore include both proofs
to further understanding of the result.

We take as our base point the $k$th roots 
of unity. Let $\gamma$ be the path from the $k$th roots of unity to the 
$(k-1)$th given by projecting to \exps{k}\ the linear homotopy from 
$(0,1/k,\ldots,(k-1)/k)$ to $(0,0,1/(k-1),\ldots,(k-2)/(k-1))$. 
Let $\alpha$ be the loop given by rotating the $k$th roots of
unity anti-clockwise through $2\pi/k$, and $\beta$ the loop based at
$\{\lambda^k=1\}$ given by taking $\gamma$ to $\{\lambda^{k-1}=1\}$, 
rotating the circle anti-clockwise through $2\pi/(k-1)$ then following
$\gamma$ back to the basepoint. Then:

\begin{theorem}
$\pi_1(\exps{k}\setminus\exps{k-2})$ has presentation
$\langle \alpha,\beta | \alpha^k=\beta^{k-1} \rangle$.
\end{theorem}

\begin{figure}[ht!]
\begin{center}\small
\leavevmode
\psfrag{a}{$\alpha$}
\psfrag{b}{$\beta$}
\psfrag{c}{$\delta$}
\psfrag{d}{$\delta\alpha$}
\includegraphics{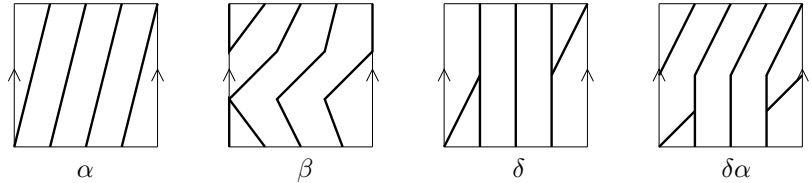}
\caption{Movies of the generators of $\pi_1(\exps{4}\setminus\exps{2})$. Each
         square represents a cylinder $\s\times[0,1]$; dark lines show
         the motion of points. The diagram for $\beta$ may be isotoped
         rel boundary (and respecting the fact that each slice $\s\times
         \{t\}$ should meet the curves in either $3$ or $4$ points) to give
         the diagram on the right, showing $\beta\simeq\delta\alpha$.}
\label{generators.fig}
\end{center}
\end{figure}

\begin{remark}
Another choice of second generator is the loop $\delta$ given by
``teleporting'' a point from $e^{2\pi i/k}$ to $e^{-2\pi i/k}$ as
follows. Between time $t=0$ and $1/2$ move one point from 
$1$ to $e^{2\pi i/k}$, keeping the rest fixed. At time $t=1/2$ the 
moving point merges with $e^{2\pi i/k}$ and we may regard it
as being at $e^{-2\pi i/k}$ instead; between $t=1/2$ and $1$ move the
extra point at $e^{-2\pi i/k}$ back to $1$, keeping the rest fixed. 
Figure~\ref{generators.fig} illustrates that $\beta$ is homotopic to 
$\delta\alpha$, and we will use this below to 
show that the inclusion
\[
\exps{k}\setminus\exps{k-2}\hookrightarrow\exps{k}\setminus\exps{k-3}
\]
is trivial on $\pi_1$.
\end{remark}

\begin{proof}
The result is an application of Van Kampen's theorem. Each of
$\exps{k}\setminus\exps{k-1}$ and $\exps{k-1}\setminus\exps{k-2}$
has the homotopy type of a circle, with fundamental group generated
(up to basepoint) by $\alpha$ and $\beta$ respectively. We show that
$\exps{k-1}\setminus\exps{k-2}$ has a neighbourhood $N$ in 
$\exps{k}\setminus\exps{k-2}$ that deformation retracts to it and
apply Van Kampen's theorem to the cover consisting of 
$\exps{k}\setminus\exps{k-1}$ and $N$. The intersection $N'$ of these
two sets will turn out to have the homotopy type of a circle also,
and will lead to the relation $\alpha^k=\beta^{k-1}$.

The preimage of $\exps{k}\setminus\exps{k-2}$ in $\Delta\times I$ is the
product of $\Delta$ less all faces of codimension two or more with the
interval. Let $b$ be the barycentre of $\Delta$, set 
\begin{eqnarray*}
\Delta_i & = & [b,v_{i+1},\ldots,v_{k-1},v_0,\ldots,v_{i-1}], \\
\delta_i & = & [v_{i+1},\ldots,v_{k-1},v_0,\ldots,v_{i-1}],
\end{eqnarray*}
and note that $\Delta$ is the union of the $\Delta_i$ and its boundary is
the union of the $\delta_i$. Further $\phi$ sends $\Delta_i,\delta_i$ to
$\Delta_{i-1},\delta_{i-1}$ respectively, where subscripts are taken 
mod $k-1$. Let
\[
N = q \biggl( \bigcup_{i=0}^{k-1} \left(\intr{\Delta_i} 
              \cup \intr \delta_i \right) \times I \biggr)
\]
and observe that $N$ is a neighbourhood of 
\[
\exps{k-1}\setminus\exps{k-2} = 
q \biggl( \bigcup_{i=0}^{k-1} \intr \delta_i  \times I \biggr)
\]
in $\exps{k}\setminus\exps{k-2}$. 

The half open simplex $\intr \Delta_i\cup\intr\delta_i$ deformation 
retracts to $\intr \delta_i$ and moreover this may be done for each
$i$ simultaneously in a way compatible with the action of $\phi$. 
Crossing this with $I$
gives a deformation retraction of 
\begin{eqnarray*}
\bigcup_{i=0}^{k-1} \left(\intr{\Delta_i} 
              \cup \intr \delta_i \right) \times I
& \mbox{to} &
\bigcup_{i=0}^{k-1} \intr \delta_i  \times I
\end{eqnarray*}
that descends to a deformation retraction of $N$ onto 
$\exps{k-1}\setminus\exps{k-2}$. Thus $\exps{k-1}\setminus\exps{k-2}$
does have a neighbourhood as desired and by the Van Kampen theorem
\[
\pi_1(\exps{k}\setminus\exps{k-2}) = 
\langle\alpha\rangle\ast_{\pi_1(N')}\langle\beta\rangle,
\]
where
\[
N'=N\cap (\exps{k}\setminus\exps{k-1}) = 
q \biggl( \bigcup_{i=0}^{k-1} \intr{\Delta_i} \times I \biggr).
\]

Now $N'$ is homeomorphic to 
$\left(\bigcup_i \intr\Delta_i\times I\right)/
\left(\Delta_i\times \{1\}\sim\Delta_{i-1}\times\{0\}\right)$,
where the gluing relations are given by $\phi$. Thus
$N'\cong \intr\Delta_{k-1}\times \s$, where we have chosen $\Delta_{k-1}$
since $q(\Delta_{k-1}\times \{0\})$ contains the path $\gamma$. If
$\pi_1(N')=\langle\varepsilon\rangle$ then clearly $\varepsilon=\alpha^k$
in $\pi_1(\exps{k}\setminus\exps{k-1})$; it remains to determine the image of
$\varepsilon$ in $\pi_1(\exps{k-1}\setminus\exps{k-2})$.

Pushing $\varepsilon$ onto $\exps{k-1}\setminus\exps{k-2}$ via the deformation 
retraction of $N$ we see that upstairs in $\Delta\times I$ $\varepsilon$
traverses the length of each $\intr\delta_i\times I$ exactly once and 
positively.
For $i=1,\ldots,k-1$ we have 
$q(\intr\delta_i\times I)=\exps{k-1}\setminus\exps{k-2}$ and we pick up a copy 
of the generator of $\pi_1(\exps{k-1}\setminus\exps{k-2})$ from each. 
However, recalling that $\delta_0\times I$ is the fake face mapping only to
\expcupone{k-1}\ we see that the contribution from this face is just the
constant loop. Thus $\varepsilon$ maps to $k-1$ times the generator in 
$\pi_1(\exps{k-1}\setminus\exps{k-2})$ and we get
the relation $\alpha^k=\beta^{k-1}$ as required.
\end{proof}

\begin{figure}[ht!]
\begin{center}\small
\leavevmode
\psfrag{a}{(a)}
\psfrag{b}{(b)}
\psfrag{c}{(c)}
\psfrag{d}{(d)}
\psfrag{e}{(e)}
\includegraphics{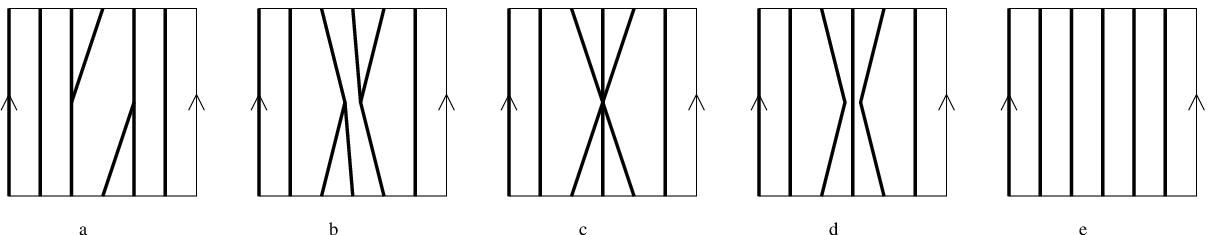}
\caption{The movie of a homotopy from $\delta$ to the constant loop
in \exps{k}\ less the codimension three strata. The figures show $k=6$
but the method clearly generalises.
Figure (a) shows $\delta$, where we have cut the circle at $-1$ for 
clarity. Planar isotopy gives (b), then we merge the branch points (c)
and separate them again (d) so that there are 6 distinct points throughout
the path. Planar isotopy leads to the constant loop in (e).
}
\label{nullhomotopic.fig}
\end{center}
\end{figure}

Both $\alpha$ and $\beta$ are null homotopic in $\exps{k}$ 
less the codimension three strata \exps{k-3}. To see this consider
figure~\ref{nullhomotopic.fig}, which shows a movie of a homotopy
in $\exps{k}\setminus\exps{k-3}$ from $\delta=\beta\alpha^{-1}$ to the 
constant loop. The relation $\alpha^k=\beta^{k-1}$ then gives
$\alpha=\beta=1$.

\subsection{The homotopy type of $\exps{k}\setminus\exps{k-2}$}

We now turn to the more delicate matter of showing that 
$\exps{k}\setminus\exps{k-2}$ deformation retracts to a subspace
having a cell structure as in figure~\ref{cellstructure.fig}. We will
of course construct the deformation retraction upstairs in $\Delta\times I$
and some care will be required to ensure it descends to the quotient.

Let $b_i$ be the barycentre of $\delta_i$. Since affine maps of
simplices take barycentres to barycentres the 
identifications~(\ref{facegluings.eq}) in section~\ref{newmodel.sec} 
glue $\{b_i\}\times I$ to
$\{b_j\}\times I$ for $i,j\not=0$ and $\phi$ glues $[b,b_i]\times\{1\}$ to 
$[b,b_{i-1}]\times\{0\}$ for each $i$, while $\{b_0\}\times I$ is of
course collapsed to $\{b_0\}\times\{0\}$. Letting $B$ be the star graph
$\bigcup_{i=0}^{k-1} [b,b_i]$ it is clear that $q(B\times I)$ may be
given a cell structure as in figure~\ref{cellstructure.fig}. 
The fact that $\exps{k}\setminus\exps{k-2}$ deformation retracts to 
this subspace is a consequence of the following technical lemma.

\begin{lemma}
There is a deformation retraction $\Phi$ from
\begin{eqnarray*}
\intr\Delta \cup \bigcup_{i=0}^{k-1} \intr\delta_i & \mbox{to} & B
\end{eqnarray*}
such that 
\begin{enumerate}
\renewcommand{\theenumi}{\alph{enumi}}
\item
$\Phi_t(\intr\delta_i)\subseteq\intr\delta_i$ for all $t$ and 
$i=0,\ldots,k-1$;
\label{starfacecondition.enum}
\item
$\Phi_t$ commutes with the action of the symmetric group $S_k$ 
on $\Delta$, where $\sigma\in S_k$ acts by the canonical map
$[v_0,\ldots,v_{k-1}]\rightarrow [v_{\sigma(0)},\ldots,v_{\sigma(k-1)}]$.
\end{enumerate}
\label{retracttostar.lem}
\end{lemma}

We will define $\Phi$ using the barycentric subdivision of $\Delta$, and
first show:

\begin{lemma}
Let $\beta_0,\ldots,\beta_n$ be affinely independent points. Then
$[\beta_0,\ldots,\beta_n]\setminus [\beta_0,\ldots,\beta_{n-2}]$ deformation
retracts to $[\beta_{n-1},\beta_n]$ via a homotopy $\Psi$ such that
\begin{equation}
\Psi_t(\delta\setminus [\beta_0,\ldots,\beta_{n-2}])
\subseteq \delta\setminus [\beta_0,\ldots,\beta_{n-2}]
\label{facecondition.eq}
\end{equation}
for all $t$ and each face $\delta$ of $[\beta_0,\ldots,\beta_n]$. 
\label{retracttoedge.lem}
\end{lemma}

\begin{proof}
The proof is by induction on $n$, the case $n=1$ being trivial. Define
a deformation retraction $\psi$ to $[\beta_1,\ldots,\beta_n]\setminus 
[\beta_1,\ldots,\beta_{n-2}]$ by
\[
\psi_t \biggl( \sum_{i=0}^n \lambda_i\beta_i \biggr) = 
(1-t)\sum_{i=0}^n \lambda_i\beta_i + 
t\sum_{i=1}^n \frac{\lambda_i}{1-\lambda_0}\beta_i \; ;
\]
this is well defined since $\lambda_0$ is never equal to $1$, and 
the co-efficient of at least one of $\beta_{n-1},\beta_n$ is nonzero
for all $t$. Moreover if $v$ is a convex combination of 
$\beta_{i_0},\ldots,\beta_{i_\ell}$ then $\psi_t(v)$ is too so the 
face condition~(\ref{facecondition.eq}) is satisfied. Applying the 
induction hypothesis to $[\beta_1,\ldots,\beta_n]\setminus 
[\beta_1,\ldots,\beta_{n-2}]$ gives the result.
\end{proof}

\begin{proof}[Proof of Lemma~\ref{retracttostar.lem}]
A typical simplex in the barycentric subdivision of $\Delta$ has the
form $[\beta_0,\ldots,\beta_{k-1}]$ where each $\beta_i$ is the barycentre
of an $i$--dimensional face of $\Delta$ containing 
$\beta_0,\ldots,\beta_{i-1}$. In particular $\beta_{k-1}=b$, $\beta_0$
is some vertex $v_i$ and $\beta_{k-2}$ is $b_j$ for some $j\not=i$. 
Deleting the codimension two faces of $\Delta$ deletes precisely
$[\beta_0,\ldots,\beta_{k-3}]$ from $[\beta_0,\ldots,\beta_{k-1}]$
and we define $\Phi$ on this simplex using the deformation retraction
$\Psi$ given by lemma~\ref{retracttoedge.lem}. 

Suppose $[\beta_0,\ldots,\beta_{k-1}]$ and $[\beta'_0,\ldots,\beta'_{k-1}]$
share a common face $[\beta_{i_0},\ldots,\beta_{i_\ell}]$. Then 
necessarily $\beta_{i_j}=\beta'_{i_j}$ for $j=1,\ldots,\ell$ and this 
face is fixed pointwise by the canonical map 
$[\beta_0,\ldots,\beta_{k-1}]\rightarrow[\beta'_0,\ldots,\beta'_{k-1}]$;
the face condition~(\ref{facecondition.eq}) in
Lemma~\ref{retracttoedge.lem} then shows that $\Phi$ is well defined. 
Condition~(\ref{starfacecondition.enum}) follows from the face 
condition~(\ref{facecondition.eq}) applied to each simplex of the form
$[\beta_0,\ldots,\beta_{k-2}]$, and the commutativity of $\Phi$ with 
the action of the symmetric group is a consequence of the fact that any 
$\sigma\in S_k$ permutes the barycentres of the $i$--dimensional faces.
\end{proof}

\begin{corollary}[implies Theorem~\ref{th.torusknot}] The space
$\exps{k}\setminus\exps{k-2}$ deformation retracts to $q(B\times I)$.
\end{corollary}

\begin{proof}
Crossing $\Phi$ with the identity gives a deformation retraction of 
\begin{eqnarray*}
\biggl(\intr\Delta \cup \bigcup_{i=0}^{k-1} \intr\delta_i\biggr)\times I 
& \mbox{to} & B\times I
\end{eqnarray*}
and we check that this is compatible with the gluings
\begin{enumerate}
\item
$\phi \co \Delta_i\times\{1\} \rightarrow \Delta_{i-1}\times\{0\}$,
\label{phi.glue}
\item
$[v_0,\ldots,\hat{v}_{i},\ldots,v_{k-1}]\times I \rightarrow
[v_0,\ldots,\hat{v}_{j},\ldots,v_{k-1}]\times I$ for $i,j\not=0$,
\label{sides.glue}
\item
$\delta_0\times I\rightarrow \delta_0\times\{0\}$.
\label{fakeface.glue}
\end{enumerate}
Compatibility with~(\ref{phi.glue}) follows from commutativity of $\Phi$ with
$S_k$; compatibility with~(\ref{sides.glue}) uses commutativity with 
$S_k$ together with condition~(\ref{starfacecondition.enum}) of 
Lemma~\ref{retracttostar.lem}; and compatibility with~(\ref{fakeface.glue})
comes from constructing the homotopy by crossing $\Phi$ with the identity 
on $I$. 
\end{proof}

\section{The degree of an induced map}

A map $f\co\s \rightarrow \s$ induces a map $\Exp{k}{f}\co\exps{k}\rightarrow
\exps{k}$ of homotopy spheres, the degree of which depends only on $k$ and the 
degree of $f$. We claim (Theorem~\ref{th.degree}) that
\[
\deg \Exp{k}{f} = \left(\deg f\right)^{\left\lfloor\frac{k+1}{2}\right\rfloor}.
\]

\begin{proof}[Proof of Theorem~\ref{th.degree}]
We begin by reducing to the case where $k=2\ell -1$ is odd, using the
commutative diagram
\[
\begin{CD}
\exps{2\ell-1} @>\Exp{2\ell-1}{f}>> \exps{2\ell-1} \\
@VV\cup\{1\}V                 @VV\cup\{f(1)\}\; \; .V\\
\exps{2\ell} @>\Exp{2\ell}{f}>> \exps{2\ell}  
\end{CD}
\]
The vertical arrows are degree one maps by the results of 
section~\ref{expkS1.sec}, so we have
\[
\deg \Exp{2\ell}{f} = \deg \Exp{2\ell-1}{f}
\]
and it suffices to show that
\mbox{$\deg \Exp{2\ell-1}{f} = \left(\deg f\right)^\ell$}.
We do this by considering separately the cases $\deg f > 0$ and 
$\deg f = -1$; in both cases we assume that $f$ is the map $\lambda\mapsto
\lambda^d$ and count the preimages of a generic point
with signs.

Suppose that $f(\lambda)=\lambda^d$ with $d$ positive. A generic point
$\Lambda$ of \exps{k}\ will have $d^k$ preimages under $\Exp{k}{f}$, 
corresponding to the $d$ choices for the preimage of each element of $\Lambda$.
For concreteness let $\lambda_0, \ldots, \lambda_{k-1}$ be the $k$th roots
of unity, cyclicly ordered so that
\[
\lambda_r = e^{2\pi ir/k}.
\]
Under $f$ each $\lambda_r$ has $d$ preimages $\lambda_{r,s}$, 
$s=0,\ldots,d-1$, which we again cyclicly order so that
\[
\lambda_{r,s} = e^{2\pi i(sk+r)/kd}. 
\]
Then a preimage of $\{\lambda|\lambda^k=1\}$ is specified by a $k$--tuple
$(s_0,\ldots,s_{k-1})$ of integers mod $d$, corresponding to the set
\[
\{\lambda_{0,s_0},\ldots,\lambda_{k-1,s_{k-1}}\},
\]
and comes with a positive or negative sign according to whether
this set is an even or odd permutation of cyclic order when ordered as
written. Note that the sign relative to cyclic ordering makes sense
since $k=2\ell -1$ is odd. 

\begin{figure}[ht!]
\begin{center}\small
\leavevmode
\psfrag{lam00}{$\lambda_{0,s_0}$}
\psfrag{lam01}{$\lambda_{0,s_1}$}
\psfrag{lam10}{$\lambda_{1,s_0}$}
\psfrag{lam11}{$\lambda_{1,s_1}$}
\includegraphics{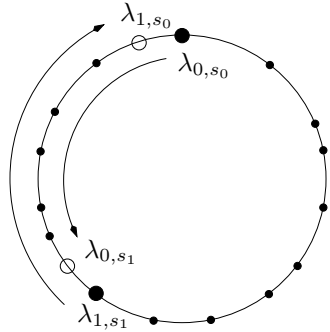}
\caption{The preimages $S$ and $\sigma_{01}(S)$, if different, have opposite 
         signs. We may move from one to the other by hopping a point at
         $\lambda_{0,s_0}$ around the circle to $\lambda_{0,s_1}$, and
         then hopping a point at $\lambda_{1,s_1}$ around the circle to 
         $\lambda_{1,s_0}$. If we move both points around the same arc then
         both move past the same points in between, except that one must
         pass the other as well. This requires an odd number of
         transpositions.}
\label{switching.fig}
\end{center}
\end{figure}

Our goal is to match the preimages
up in cancelling pairs until those that are left are all positive. To
this end consider the involution $\sigma_{01}$ given by
$(s_0,\ldots,s_{k-1})\mapsto (s_1,s_0,s_2,\ldots,s_{k-1})$; we claim 
that if a preimage $S$ is not fixed by $\sigma_{01}$ then $S$ and 
$\sigma_{01}(S)$ have opposite signs. Indeed, if $s_0\not=s_1$ and
we move points from 
$\lambda_{0,s_0}$ to $\lambda_{0,s_1}$ and from $\lambda_{1,s_1}$ to 
$\lambda_{1,s_0}$
around the same arc of the circle then (since there can be no points between
$\lambda_{0,s_j}$ and $\lambda_{1,s_j}$) both must pass exactly the same
points in between, and in addition one must pass the other 
(see figure~\ref{switching.fig}). 
This involves an odd number of transpositions
so $S$ and $\sigma_{01}(S)$ have opposite signs if $s_0\not=s_1$.

Applying the same argument in turn to the involutions switching 
$s_{2j}$ and  $s_{2j+1}$, $1\leq j\leq \ell-2$,
acting just on those 
preimages fixed by all previous involutions, we see that we may match up all 
preimages in cancelling pairs except those for which $s_{2j}=s_{2j+1}$,
$0\leq j\leq \ell-2$. Since these can all be shuffled to cyclic order
by moving points around in pairs they are all positive, and there are
$d^\ell$ of them as there are $d$ choices for each $s_{2j}$, 
$j=0,\dots,\ell-1$. This gives $\deg \Exp{k}{f}=d^\ell$ as desired.

Now consider $f(\lambda)=\lambda^{-1}$. The sole preimage of the
ordered set $\{\lambda_0,\ldots,\lambda_{k-1}\}$ is
$\{\lambda_0,\lambda_{k-1},\ldots,\lambda_1\}$
which may be put in cyclic order using $(k-1)/2=\ell-1$ transpositions. 
We get an additional factor of $(-1)^k$ from the product of the local degrees 
of $f$ at each $\lambda_j$, so 
\[
\deg \Exp{k}{f} = (-1)^{\ell-1}\cdot(-1)^{2\ell-1}
= (-1)^{\ell-2} = \left(\deg f\right)^\ell 
\]
as required.

If $\deg f=0$ then clearly $\deg \Exp{k}{f}=0$, so putting the cases 
$\deg f > 0$ and $\deg f = -1$ together using 
$\deg (g\circ h) = (\deg g)\cdot(\deg h)$ gives the result.
\end{proof}

\appendix

\section{Cut and paste proof of Theorems~\ref{th.Bott} and~\ref{th.trefoil}}
\label{cutandpaste.apdx}

\begin{figure}[ht!]
\begin{center}\small
\leavevmode
\psfrag{key:}{key:}
\psfrag{A}{$A$}
\psfrag{B}{$B$}
\psfrag{C}{$C$}
\psfrag{D}{$D$}
\psfrag{E}{$E$}
\psfrag{F}{$F$}
\psfraga <0pt, -2pt> {1}{visible edges of handlebody}
\psfraga <0pt, -2pt> {2}{hidden edges of handlebody}
\psfraga <0pt, -2pt> {3}{visible segments of $\alpha$}
\psfraga <0pt, -2pt> {4}{hidden segments of $\alpha$}
\psfraga <0pt, -2pt> {5}{back segment of $\alpha$}
\psfraga <0pt, -2pt> {6}{attaching curve $\beta$}
\includegraphics[width=3.5in]{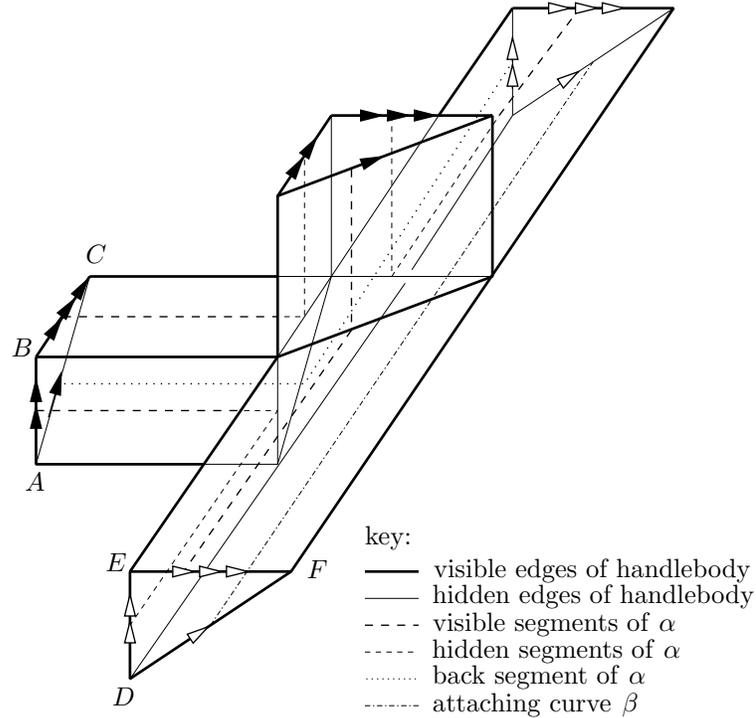}
\caption{A regular neigbourhood of the dual $1$--skeleton of the 
triangulation of \exps{3}\ in figure~\ref{3simplex.fig}. Glue the
triangular faces as indicated to get a genus two solid handlebody
$H$ forming half of a Heegaard splitting of \exps{3}. The curves
$\alpha$ and $\beta$ bound discs in the second handlebody $H'$.}
\label{heegaard.fig}
\end{center}
\end{figure}

We give a mostly pictorial proof that \exps{3}\ is $S^3$ and that
\exps{1}\ inside it is a left-handed trefoil knot.

\begin{figure}[ht!]
\begin{center}\small
\leavevmode
\psfrag{a (hidden)}{$\alpha$ (hidden)}
\psfrag{a (visible)}{$\alpha$ (visible)}
\psfrag{b}{$\beta$}
\psfraga <-2pt, 0pt> {A}{$A$}
\psfrag{B}{$B$}
\psfrag{C}{$C$}
\psfrag{D}{$D$}
\psfrag{E}{$E$}
\psfrag{F}{$F$}
\includegraphics[width=2.5in]{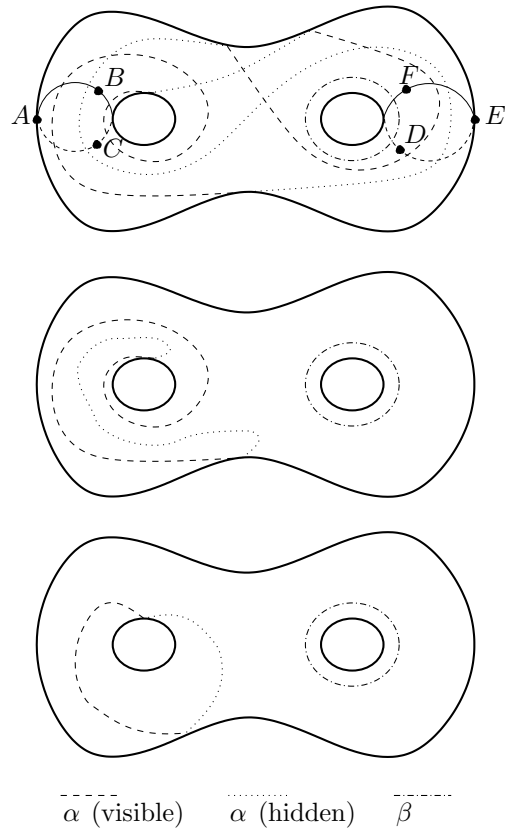}
\caption{Heegaard diagrams of \exps{3}.
         Top: The handlebody $H$ obtained by gluing the
         triangular faces of figure~\ref{heegaard.fig}. 
         Middle: Attaching curve $\alpha$ after sliding 
         the loops going around the right $1$--handle over the cancelling
         $2$--handle. Bottom: After isotopy $\alpha$ forms 
         a $(1,1)$--curve around the left $1$--handle, giving a Heegaard
         diagram of $S^3$.}
\label{genus2.fig}
\end{center}
\end{figure}

We found in section~\ref{exp3.sec} that \exps{3}\ is a $3$--manifold with
a triangulation consisting of just one $3$--simplex, and in the standard way
we obtain a Heegaard splitting by regarding it as the union 
$H'\cup_\partial H$ of regular neighbourhoods of the $1$-- and dual 
$1$--skeletons. A regular neighbourhood of the dual $1$--skeleton of the
$3$--simplex of figure~\ref{3simplex.fig}(a) is shown in 
figure~\ref{heegaard.fig} and
$H$ is the genus two solid handlebody given by gluing the four triangular faces
at the end of each ``arm'' as indicated by the arrows. We keep track
of $H'$ by recording loops $\alpha$ and $\beta$ linking each of the edges,
the loops forming the attaching circles for the $2$--handles of \exps{3}. The
loop $\alpha$ linking the edge $a$ of figure~\ref{3simplex.fig}(b) appears
in five pieces in figure~\ref{heegaard.fig}, four of which are
shown dashed and the fifth dotted, while the loop $\beta$ linking the edge
$b$ appears in just one piece, indicated by the dash-dot-dash-dot line.

\begin{figure}[ht!]
\begin{center}\small
\leavevmode
\psfrag{A}{$A$}
\psfrag{B}{$B$}
\psfrag{C}{$C$}
\psfrag{D}{$D$}
\psfrag{E}{$E$}
\psfrag{F}{$F$}
\includegraphics[width=2.85in]{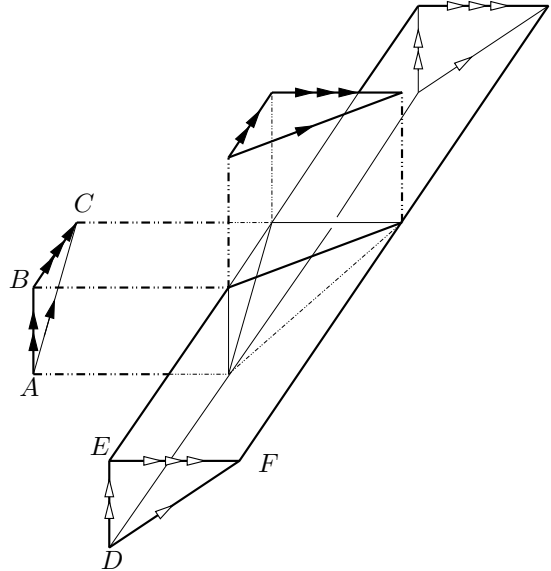}
\caption{Generic orbits of the \s\ action on \exps{3}\ meet the 
$3$--simplex in three lines parallel to the direction $(1,1,1)$. 
Perturbing them to lie on the Heegaard surface we may obtain the
three arcs shown dash-dot-dotted.}
\label{cutorbit.fig}
\end{center}
\end{figure}

Bending the arms and giving the top arm a one-third twist so the arrows
match we glue the triangular faces to obtain $H$, shown in the top 
diagram in figure~\ref{genus2.fig}. 
The curve $\beta$ is again
indicated by a dash-dot-dash-dot pattern and we see immediately that the 
attached $2$--handle cancels the right $1$--handle. The curve $\alpha$ 
linking $a$, shown dashed when it is on top and dotted when it is underneath,
goes geometrically twice over the right $1$--handle and we slide each loop
over the cancelling $2$--handle to get the middle diagram. Further isotopy
leads to the bottom diagram in which $\alpha$ appears as a $(1,1)$--curve
around the left $1$--handle. This may be recognised as a Heegaard diagram
for $S^3$ but we nevertheless apply a Dehn twist to convert it to 
the standard diagram. Writing the meridian first and giving $\partial H$
the right-hand orientation induced by $H$ the appropriate Dehn twist
acts by 
\[
T=\twobytwo{1}{-1}{0}{1}
\]
on $H_1(\partial(H\cup_\beta (D^2\times I)))$.

\begin{figure}[ht!]
\begin{center}\small
\leavevmode
\psfrag{A}{$A$}
\psfrag{B}{$B$}
\psfrag{C}{$C$}
\psfrag{D}{$D$}
\psfrag{E}{$E$}
\psfrag{F}{$F$}
\includegraphics[width=2.8in]{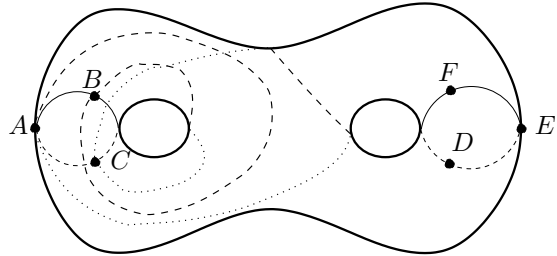}
\caption{The perturbed orbit, shown dashed when on top and dotted when
below, after gluing figure~\ref{cutorbit.fig} up to form $H$. It traces
a $(1,3)$--curve on the left $1$--handle and is transformed
to a $(-2,3)$--curve by the Dehn twist.}
\label{orbit.fig}
\end{center}
\end{figure}

We now turn our attention to $\exps{1}\hookrightarrow\exps{3}$,
corresponding to the edge $b$ of figure~\ref{3simplex.fig}(b). We take
a push-off of $b$ into the interior of the simplex and perturb it to lie
on the Heegaard surface. This is most easily done by recalling that
\exps{1}\ forms a generic (here meaning trivial stabiliser) orbit of the \s\
action on \exps{3}, and that generic orbits passing through the interior
of the simplex break into three lines parallel to the vector $(1,1,1)$. 
Pushing them onto the Heegaard surface we obtain the arcs shown
dash-dot-dotted in figure~\ref{cutorbit.fig}. 
The resulting curve on 
$\partial H$ appears in figure~\ref{orbit.fig}, shown dashed when it is on top 
and dotted underneath, and as might be expected from figure~\ref{cutorbit.fig}
it forms a $(1,3)$--curve around the left $1$--handle. The Dehn twist
with matrix $T$ that converts the bottom Heegaard diagram of 
figure~\ref{genus2.fig} to the standard diagram takes this curve
to a $(-2,3)$--curve, giving a left-handed trefoil as promised.

\Addresses\recd

\end{document}